\documentclass[a4paper]{amsart}
\usepackage{amsthm,amsmath,amssymb,mathrsfs,graphicx, verbatim, overpic}
\usepackage[latin1]{inputenc}
\usepackage{tabmac}
\usepackage{graphicx, color, overpic}
\usepackage[usenames,dvipsnames,svgnames,table]{xcolor}

\usepackage{tikz}

 \newtheorem{thm}{Theorem}

 \theoremstyle{definition}
 
 \newtheorem{ex}[thm]{Example}

 \theoremstyle{remark}
 
\numberwithin{thm}{section}

\def\R{{\mathbb R}}
\def\Z{{\mathbb Z}}

\def\C{{\mathbb C}}

\newcommand{\liz}{\textcolor{blue}}
\newcommand{\tr}{\textcolor{red}}

\newcommand{\height}{\operatorname{ht}}

\definecolor{ltgreen}{rgb}{0.0, 0.5, 0.0}
\definecolor{dkgreen}{rgb}{0.0, 0.42, 0.24}

\begin{document}

\begin{title}
{Reduced words for reflections in Weyl groups}
\end{title}
\author{Elizabeth Mili\'{c}evi\'{c}}

\address{Elizabeth Mili\'{c}evi\'{c}, Department of Mathematics \& Statistics, Haverford College, 370 Lancaster Avenue, Haverford, PA, 19041, USA}
\email{emilicevic@haverford.edu}

\thanks{The author was partially supported by NSF Grants DMS-1600982 and DMS-2202017.}

\begin{abstract}
The reflections in a Coxeter group are defined as conjugates of a single generator, and thus admit palindromic expressions as products of generators. Our main result gives closed formulas providing a palindromic reduced expression for each reflection in any finite Weyl group. There exist algorithmic methods for determining such reduced expressions, but explicit formulas have not been recorded outside of well-known special cases.
\end{abstract}

\maketitle

\section{Introduction}\label{intro}

Given any Coxeter group $W$ with set of generators $S$, every element $w \in W$ can be written as a product $w=s_{1}s_2\cdots s_{k}$ for some $s_{i} \in S$.  If $k$ is minimal, this expression for $w$ is reduced.  Every $w \in W$ admits at least one reduced expression, and all expressions can be transformed into a reduced expression by doing a sequence of nil, braid, and/or commuting relations on the generators. Identifying reduced expressions with desirable properties, as well as enumerating reduced expressions, have motivated substantial work on the combinatorics of Coxeter groups.

\subsection{Reduced words for reflections}

Every Coxeter group has a distinguished subset of reflections $T = \{ wsw^{-1} \mid s \in S, w \in W\}$. Equivalently, every $t \in T$ has a palindromic expression in the generators, meaning that the product reads the same forwards or backwards.  It is a straightforward exercise to prove that every reflection admits a palindromic expression which is also reduced; see Exercise 1.10 in \cite{BB} and a proof in either Lemma 4.1 of \cite{BFP} or the response to this question on MathOverflow \cite{Samuel}.  Explicit formulas for reduced expressions for reflections do not appear in the literature on Coxeter groups, however, apart from well-known special cases such as $W = S_{n}$ the symmetric group.

We now restrict to the class of Coxeter groups which arise as Weyl groups of a crystallographic root system.  This family plays a central role in the classification of semisimple Lie algebras, and almost all finite irreducible Coxeter groups are Weyl groups.   If $W$ is a Weyl group, each $t \in T$ acts geometrically on Euclidean space as a reflection in the hyperplane orthogonal the corresponding root vector. The purpose of this paper is to provide closed formulas for palindromic reduced expressions for all reflections in any Weyl group.

\begin{thm}\label{thm:main}The following are palindromic reduced expressions, one for each distinct reflection in the Weyl group of the type specified.
\begin{enumerate}
\item[($A_n$)]  For $1 \leq i < j \leq n$, we have
$$s_{e_i-e_{j+1}} = s_is_{i+1} \cdots s_{j-1}s_j s_{j-1} \cdots s_{i+1}s_i. \phantom{xxxxxxxxxxxxxxxxxxxxxxxxxxxxxxxxx}$$
\item[($B_n$)] For $1 \leq i < j \leq n$, we have 
\begin{align*} s_{e_i-e_j} &= s_is_{i+1} \cdots s_{j-2}s_{j-1} s_{j-2} \cdots s_{i+1}s_i \\
 s_{e_i+e_j} &= \left(s_js_{j+1} \cdots s_{n-1} \cdot s_is_{i+1} \cdots s_{n-2}\right)\left(s_n s_{n-1} s_n\right) \left( s_{n-2} \cdots s_{i+1}s_i \cdot s_{n-1} \cdots s_{j+1} s_j\right)\\
  s_{e_{i}} & = s_is_{i+1} \cdots s_{n-1}s_n s_{n-1} \cdots s_{i+1}s_i.
\end{align*}
\item[($C_n$)] For $1 \leq i < j \leq n$, we have 
\begin{align*} s_{e_i-e_j} &= s_is_{i+1} \cdots s_{j-2}s_{j-1} s_{j-2} \cdots s_{i+1}s_i \\
 s_{e_i+e_j} &= \left(s_js_{j+1} \cdots s_{n-1} \cdot s_is_{i+1} \cdots s_{n-2}\right)\left(s_n s_{n-1} s_n\right) \left( s_{n-2} \cdots s_{i+1}s_i \cdot s_{n-1} \cdots s_{j+1} s_j\right)\\
  s_{2e_{i}} & = s_is_{i+1} \cdots s_{n-1}s_n s_{n-1} \cdots s_{i+1}s_i.
\end{align*}

\item[($D_n$)] We have 
\begin{align*} s_{e_i-e_j} &= s_is_{i+1} \cdots s_{j-2}s_{j-1} s_{j-2} \cdots s_{i+1}s_i \ \  \text{for}\ 1 \leq i < j \leq n, \\
 s_{e_i+e_n} &=  s_is_{i+1} \cdots s_{n-2}s_n s_{n-2} \cdots s_{i+1}s_i \ \ \text{for}\ 1 \leq i < n,\\
s_{e_i+e_j} &= \left(s_js_{j+1} \cdots s_{n-2} \cdot s_i s_{i+1}\cdots s_{n-3}\right)\left(s_{n-1}s_{n-2}s_ns_{n-2}s_{n-1}\right) \left( s_{n-3} \cdots s_i \cdot s_{n-2} \cdots s_j\right)\\
 & \phantom{=ab} \text{for} \ 1 \leq i < j < n.
\end{align*}
\item[($E_6$)] See Table \ref{tab:E6}.
\item[($E_7$)] See Tables \ref{tab:E6} and \ref{tab:E7}.
\item[($E_8$)] See Tables \ref{tab:E6}, \ref{tab:E7}, \ref{tab:E8-1}, and \ref{tab:E8-2}.
\item[($F_4$)] See Table \ref{tab:F4}.
\item[($G_2$)] See Table \ref{tab:G2}.
\end{enumerate}
\end{thm}
\noindent The expansion for each root vector in terms of a basis of simple roots is also provided in Section \ref{sec:classical} for the classical types, and in each table in Section \ref{sec:exceptional} for the exceptional types.  It is the hope of the author that these formulas will be viewed as an additional reference expanding upon the helpful plates at the end of \cite{Bour46}. We remark that we also correct 5 distinct typos in item (II) of those plates within the body of this work, in each of types $C_n, D_n, E_8,$ and $G_2$.

There exist algorithmic methods for determining reduced expressions for general $w \in W$, such as the numbers game; see Section 4.3 of \cite{BB}.  In response to a Stack Exchange question explicitly requesting the results of Theorem \ref{thm:main}, Reading explains a method for starting with a word for $wsw^{-1}$ and reducing it while maintaining its palindromic property.  He also provides an alternate algorithm using the geometric realization of $(W,S)$; see \cite{Reading}.  Lemma 1.7 in \cite{BonDyer} provides a closed formula for a palindromic reduced expression of $t \in T$, given another reduced expression for $t$ as input.  Several computational software packages have implemented similar algorithms, such as the \verb!associated_reflection()! command in Sage, which returns a reduced word given a linear combination of simple roots \cite{sage}. The Maple packages \verb!coxeter! and \verb!weyl! by Stembridge support the \verb!vec2fc! command, which performs the geometric algorithm in \cite{Reading}.

\subsection{Discussion of the results}

We continue by illustrating Theorem \ref{thm:main} in types $C_4$ and $E_8$. 

\begin{ex}\label{ex:introC4}
Let $e_i$ denote the $i^{\text{th}}$ standard basis vector in $\R^4$, and choose the simple roots to be $\Delta = \{\alpha_1,\alpha_2,\alpha_3,\alpha_4\} =  \{e_1-e_2, \; e_2-e_3, \; e_3-e_4, \; 2e_4\}$. The set $\Delta$ corresponds to the simple reflections $\{s_1, s_2, s_3, s_4\}$, labeling the Dynkin diagram as follows:  
 \begin{center}
 \begin{tikzpicture}[scale=0.7,double distance=3pt,thick]
\begin{scope}[scale=1.2] 
\fill (0,-.1) circle (3pt); 
\fill (1,-.1) circle (3pt); 
\fill (2,-.1) circle (3pt); 
\fill (3,-.1) circle (3pt); 
\draw (0,-.1)--(1,-.1)--(2,-.1);
\draw (2,-0.05)--(3,-0.05);
\draw (2,-0.15)--(3,-0.15);
\draw (2.6,0.05)--(2.45,-.1);
\draw (2.45,-.1)--(2.6,-0.25);
\draw (0,-0.4) node {$s_1$};	
\draw (1,-0.4) node {$s_2$};	
\draw (2,-0.4) node {$s_{3}$};	
\draw (3,-0.4) node {$s_{4}$};	
\end{scope}
\end{tikzpicture}
\end{center}

The $4^2 = 16$ positive roots with respect to $\Delta$ are 
 \[\Phi^+=\{e_i-e_{j}\mid 1 \leq i<j \leq 4\} \cup \{ e_i+e_j \mid 1 \leq i < j \leq 4 \} \cup \{ 2e_i \mid i \in [4]\}.\]  
 The $4 \choose 2$ $=6$ positive roots in the first set are the roots for the underlying type $A_{3}$ root system generated by the first $3$ nodes on the Dynkin diagram.  Those non-simple type $A_3$ reflections thus have the following familiar palindromic reduced expressions
\[ s_{e_1-e_3} = s_1s_2s_1, \quad  s_{e_2-e_4} = s_2s_3s_2, \quad  s_{e_1-e_4} = s_1s_2s_3s_2s_1. \]
There are 3 additional non-simple reflections of type $A_4$, all conjugates of the generator $s_4$ to avoid duplication:
\[ s_{2e_{1}} = s_1s_{2} s_3s_4 s_{3}s_{2}s_1, \quad s_{2e_{2}} = s_{2} s_3s_4 s_{3}s_{2}, \quad s_{2e_{3}} = s_3s_4 s_{3}. \]
 
 The remaining 6 reflections are more interesting, and also fully illustrate the key idea for recording palindromic reduced expressions for all the non-type $A$ reflections.  First note that there is one obvious unused palindromic reduced expression $s_4s_3s_4$, which is not equivalent to $s_3s_4s_3$ in type $C_4$.  
 The remaining reflections are all conjugates of $s_4s_3s_4$ by the minimal length coset representatives of the quotient of the subgroup $\langle s_1, s_2, s_3\rangle$ by the maximal parabolic subgroup $\langle s_1,s_3\rangle$.  These representatives index the Schubert cells in the Grassmannian $Gr(2,4)$ of 2-planes in $\C^4$, and are in bijection with the Young diagrams contained in a $2 \times 2$ rectangle. A corresponding reduced word is read by overlaying the Young diagram on the grid labeled by
\[
  \tableau[Ys]{ \liz{2}& \tr{3} \\ \liz{1} & \tr{2}\\ } \quad \longleftrightarrow \quad s_{\tr{2}}s_{\tr{3}}s_{\liz{1}}s_{\liz{2}}.
\]
We build a word by reading the labels up each column, starting with the rightmost column, as shown in the figure above.  The root $e_i+e_j$ corresponds to the diagram whose first column extends down to label $i$ and second column to $j$. For example, the rectangle above corresponds to the root $e_1+e_2$, and the associated reflection is $s_{e_1+e_2} = s_2s_3s_1s_2\cdot s_4s_3s_4 \cdot s_2s_1s_3s_2$.  Likewise, 
\[ s_{e_3+e_4} = s_4s_3s_4, \quad  s_{e_2+e_4} = s_2\cdot s_4s_3s_4 \cdot s_2, \quad   s_{e_1+e_4} = s_1s_2\cdot s_4s_3s_4 \cdot s_2s_1, \]
\[ s_{e_2+e_3} = s_3s_2 \cdot s_4s_3s_4 \cdot s_2s_3, \quad   s_{e_1+e_3} = s_3s_1s_2\cdot s_4s_3s_4 \cdot s_2s_1s_3, \]
which completes our list of palindromic reduced expressions for the 16 reflections in type $C_4$.
\end{ex}

    \begin{figure}[h]
\begin{center}
 \resizebox{3.3in}{!}
{
\begin{overpic}{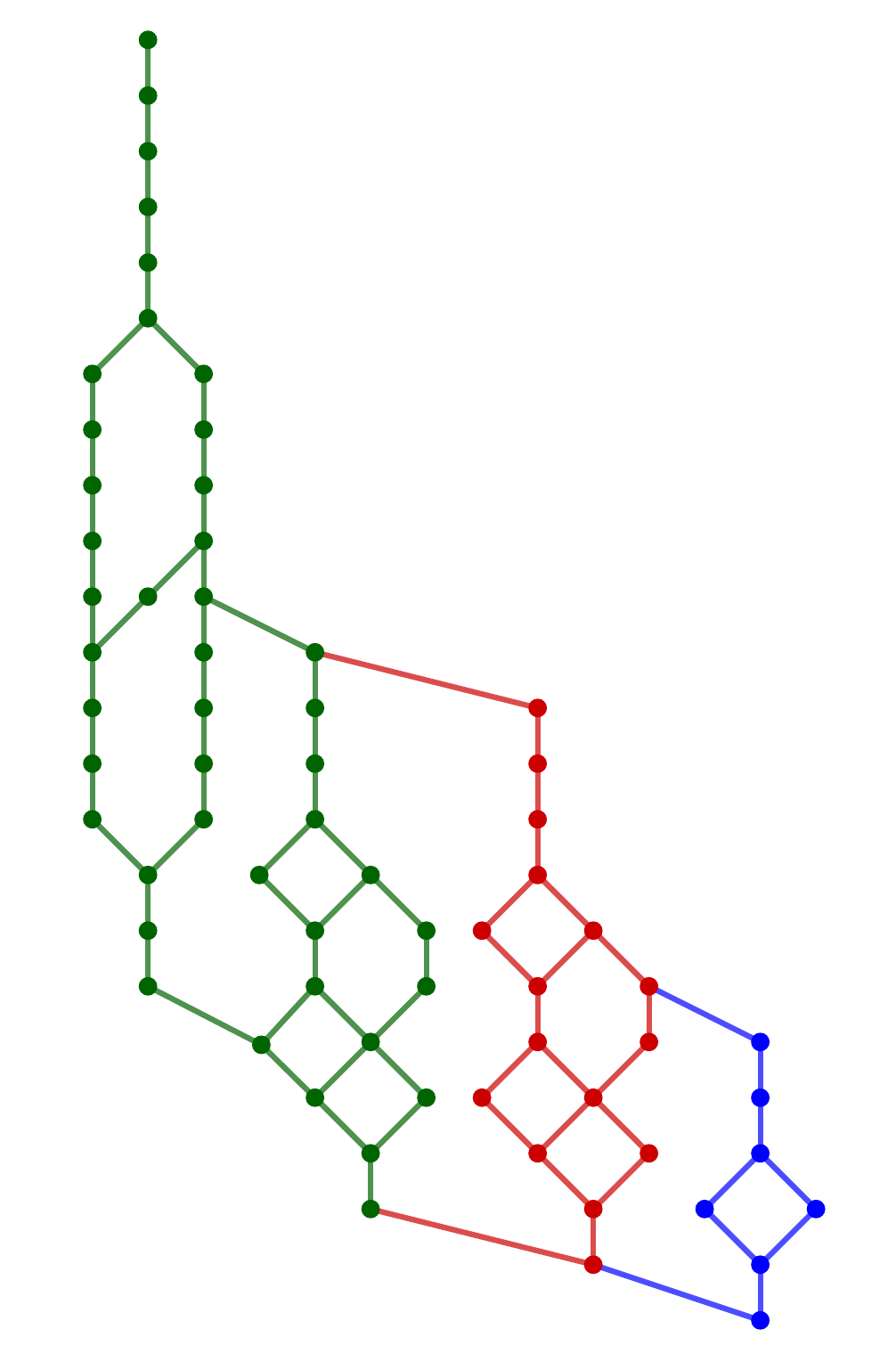}
\put(9.75,98.75){\bf \huge{$s_{\tilde{\alpha}}$}}
\put(9,94.5){\large{8}}
\put(9,90.5){\large{7}}
\put(9,86.5){\large{6}}
\put(9,82.5){\large{5}}
\put(9,78.5){\large{4}}
\put(7.5,75){\large{3}}
\put(13.5,75){\large{2}}
\put(5,70.25){\large{1}}
\put(16,70.25){\large{3}}
\put(5,66.25){\large{2}}
\put(16,66.25){\large{4}}
\put(5,62.25){\large{4}}
\put(16,62.25){\large{5}}
\put(5,58.25){\large{3}}
\put(16,58.25){\large{6}}
\put(11.25,58.25){\large{1}}
\put(5,54.05){\large{5}}
\put(16,54.05){\large{1}}
\put(10,54.05){\large{3}}
\put(5,50.1){\large{4}}
\put(16,50.1){\large{3}}
\put(5,46){\large{2}}
\put(16,46){\large{4}}
\put(5,41.85){\large{6}}
\put(16,41.85){\large{5}}
\put(7.5,37){\large{5}}
\put(13.5,37){\large{2}}
\put(9,33.75){\large{4}}
\put(9,29.75){\large{3}}
\put(21,54.05){\large{7}}
\put(21,50.1){\large{1}}
\put(21,46){\large{3}}
\put(21,41.85){\large{4}}
\put(19.5,38.25){\large{2}}
\put(25.75,38.25){\large{5}}
\put(19.5,33.25){\large{5}}
\put(25.75,33.35){\large{2}}
\put(30.05,33.45){\large{6}}
\put(21,29.85){\large{4}}
\put(15.5,26.05){\large{7}}
\put(19.5,26.05){\large{3}}
\put(25.75,25.9){\large{6}}
\put(19.5,21.05){\large{6}}
\put(25.75,21.15){\large{3}}
\put(29.35,29.85){\large{2}}
\put(30.35,25.9){\large{4}}
\put(30.35,21.25){\large{5}}
\put(30.25,17.15){\large{3}}
\put(23.45,17){\large{5}}
\put(25.75,13.15){\large{4}}
\put(30.35,51.15){\large{8}}
\put(37.25,46.1){\large{1}}
\put(37.25,42){\large{3}}
\put(37.25,37.85){\large{4}}
\put(35.75,34.25){\large{2}}
\put(42,34.25){\large{5}}
\put(35.75,29.25){\large{5}}
\put(42,29.35){\large{2}}
\put(46.35,29.45){\large{6}}
\put(37.25,25.85){\large{4}}
\put(35.75,22.05){\large{3}}
\put(42,21.9){\large{6}}
\put(35.75,17.05){\large{6}}
\put(42,17.15){\large{3}}
\put(45.6,25.85){\large{2}}
\put(46.6,21.85){\large{4}}
\put(46.6,17.25){\large{5}}
\put(46.5,13.15){\large{3}}
\put(39.7,13){\large{5}}
\put(42,9.15){\large{4}}
\put(34.5,10.5){\large{8}}
\put(23,9.25){\bf \huge{$s_{87}s_{\theta}s_{78}$}}
\put(51.5,26.25){\large{7}}
\put(56.5,21.5){\large{2}}
\put(56.5,17.5){\large{4}}
\put(52,13.7){\large{3}}
\put(58.25,13.55){\large{5}}
\put(52,8.7){\large{5}}
\put(58.25,8.8){\large{3}}
\put(56.5,5.25){\large{4}}
\put(48.5,6.75){\large{7}}
\put(55,1.25){\bf \huge{$s_{\theta}$}}
\end{overpic}
}
 \caption{Some conjugacy relations among type $E$ reflections, where $s_{\theta} = s_{16524342561}$.}\label{fig:conj}
\end{center}
\end{figure}

In the exceptional types, there are always embedded copies of Dynkin 
diagrams for lower-rank classical types, permitting use of those results to collect palindromic reduced 
expressions for many reflections, perhaps after relabeling.  Each exceptional type also has reflections 
distinct to that type, however, which must be identified systematically in some other way.  

\begin{ex}\label{ex:introE8}
Figure \ref{fig:conj} on the next page illustrates our method for identifying those reflections which are distinct to type $E$.   The vertices of the graph in Figure \ref{fig:conj} correspond to the 67 reflections distinct to types $E_6, E_7,$ and $E_8$, colored blue, red, and green, respectively, with length increasing going upward. The type $E_6$ reflection of minimal length is  $s_{\theta} = s_{16524342561}$, where we write $s_is_js_k = s_{ijk}$ for brevity.  

The edge labels in Figure \ref{fig:conj} indicate the simple reflection with which to conjugate to obtain the reflection indexing the adjacent vertex. 
 For example, the lowest red vertex represents the element $s_7 s_\theta s_7$.  By traveling along the lefthand side of Figure \ref{fig:conj}, for example, we also see that
 \[ s_{\tilde{\alpha}} = s_{87654312435426543765487} \cdot s_{16524342561} \cdot s_{78456734562453421345678}, \]
 where $\tilde{\alpha}$ is the highest root in type $E_8$.   Not all possible conjugacy relations are depicted in Figure \ref{fig:conj}; however, a sufficient subset of relations is provided to demonstrate a clear path from $s_\theta$ to each of the other 66 reflections distinct to type $E$.
\end{ex}

\subsection*{Acknowledgments.} Theorem \ref{thm:main} was obtained in service of a preliminary case-by-case proof of Lemma 4.9 in \cite{MilVie}, and has thus benefitted from feedback by Eva Viehmann. The author is grateful to Anne Thomas, whose tikz code was adapted to create the Dynkin diagrams displayed in this paper. Part of this work was carried out while the author was a Director's Mathematician in Residence at the Budapest Semesters in Mathematics, jointly supported by the R\'enyi Alfr\'ed Matematikai Kutat\'oint\'ezet. Calculations in some exceptional types were 
assisted by Sage, and the author thanks the Sage developers for implementing related procedures \cite{sage}.

\section{Background} \label{sec2}

In this section, we review all preliminary definitions, notation, and results we need concerning Coxeter groups and their associated root systems.  For more detailed references, we refer the reader to textbooks such as \cite{BB,Bour46,Humph}.

\subsection{Coxeter groups}

Let $(W,S)$ be a finite irreducible Coxeter system of rank $n$. The generating set $S = \{s_1, \dots, s_n\}$ are the simple reflections, and we denote the index set by $[n] = \{1, \dots, n\}$. Further suppose that $W$ is the Weyl group of a crystallographic root system; equivalently, let $W$ be a finite Weyl group.  The Dynkin diagram for $(W,S)$ can then be classified as one of the families displayed in Table 1 on the next page.

Every element $w \in W$ can be written as a product $w=s_{i_1}s_{i_2}\cdots s_{i_k}$ for some $s_{i_j} \in S$.  If $k$ is minimal, then this product for $w$ is called a reduced expression. In case of a reduced expression, we say that the length of $w$ equals $k$, denoted $\ell(w)=k$. 

 Denote by $T = \{ wsw^{-1} \mid w \in W, s \in S\}$ the set of reflections of $(W,S)$.  By choosing any expression for $w=s_{i_1}\cdots s_{i_k}$, we can decompose  $wsw^{-1} =(s_{i_1}\cdots s_{i_k})s(s_{i_k} \cdots s_{i_1})$ into a palindromic expression, meaning that the sequence of generators reads the same forwards and backwards.  Note that such a palindromic expression is not guaranteed to be reduced.

\begin{table}[ht]
\begin{center} 
\begin{tabular}{|c|c||c|c|} 
\hline
Type & Dynkin Diagram & Type & Dynkin Diagram  \\ 
\hline \hline
$A_n,\ n \geq 1$ 
&
\begin{tikzpicture}[scale=0.7,double distance=3pt,thick]
\begin{scope}[scale=1.2] 
\fill (0,.1) circle (3pt); 
\fill (1,.1) circle (3pt); 
\fill (3,.1) circle (3pt); 
\fill (4,.1) circle (3pt); 
\draw (0,.1)--(1,.1)--(1.5,.1);
\draw [dashed] (1.5,.1)--(2.5,.1);
\draw (2.5,.1)--(3,.1)--(4,.1);
\draw (0,-0.3) node {$s_1$};	
\draw (1,-0.3) node {$s_2$};	
\draw (3,-0.3) node {$s_{n-1}$};	
\draw (4,-0.3) node {$s_n$};		
\end{scope}
\end{tikzpicture}
&
$E_6$ 
&
\begin{tikzpicture}[scale=0.7,double distance=3pt,thick]
\begin{scope}[scale=1.2] 
\fill (0,-.5) circle (3pt); 
\fill (1,-.5) circle (3pt); 
\fill (2,-.5) circle (3pt); 
\fill (3,-.5) circle (3pt); 
\fill (4,-.5) circle (3pt); 
\fill (2,-1.5) circle (3pt); 
\draw (0,-.5)--(1,-.5)--(2,-.5)--(3,-.5)--(4,-.5);
\draw (2,-.5)--(2,-1.5);
\draw (0,-0.25) node {$s_1$};	
\draw (1,-0.25) node {$s_3$};
\draw (2,-0.25) node {$s_4$};
\draw (3,-0.25) node {$s_5$};
\draw (4,-0.25) node {$s_6$};
\draw (2.35,-1.4) node {$s_2$};	
\end{scope}
\end{tikzpicture}
 \\
 \hline
$B_n,\ n \geq 2$ 
&
\begin{tikzpicture}[scale=0.7,double distance=3pt,thick]
\begin{scope}[scale=1.2] 
\fill (0,0) circle (3pt); 
\fill (1,0) circle (3pt); 
\fill (3,0) circle (3pt); 
\fill (4,0) circle (3pt); 
\fill (5,0) circle (3pt); 
\draw (0,0)--(1,0)--(1.5,0);
\draw [dashed] (1.5,0)--(2.5,0);
\draw (2.5,0)--(3,0)--(4,0);	
\draw (4,0.05)--(5,0.05);
\draw (4,-0.05)--(5,-0.05);
\draw (4.4,0.15)--(4.55,0);
\draw (4.55,0)--(4.4,-0.15);
\draw (0,-0.3) node {$s_1$};	
\draw (1,-0.3) node {$s_2$};	
\draw (3,-0.3) node {$s_{n-2}$};	
\draw (4,-0.3) node {$s_{n-1}$};	
\draw (5,-0.3) node {$s_n$};
\end{scope}
\end{tikzpicture}
&
$E_7$ 
&
\begin{tikzpicture}[scale=0.7,double distance=3pt,thick]
\begin{scope}[scale=1.2] 
\fill (0,-.5) circle (3pt); 
\fill (1,-.5) circle (3pt); 
\fill (2,-.5) circle (3pt); 
\fill (3,-.5) circle (3pt); 
\fill (4,-.5) circle (3pt); 
\fill (5,-.5) circle (3pt); 
\fill (2,-1.5) circle (3pt); 
\draw (0,-.5)--(1,-.5)--(2,-.5)--(3,-.5)--(4,-.5)--(5,-.5);
\draw (2,-.5)--(2,-1.5);
\draw (0,-0.25) node {$s_1$};	
\draw (1,-0.25) node {$s_3$};
\draw (2,-0.25) node {$s_4$};
\draw (3,-0.25) node {$s_5$};
\draw (4,-0.25) node {$s_6$};
\draw (5,-0.25) node {$s_7$};
\draw (2.35,-1.4) node {$s_2$};	
\end{scope}
\end{tikzpicture}
 \\
\hline
$C_n,\ n \geq 2$ 
&
\begin{tikzpicture}[scale=0.7,double distance=3pt,thick]
\begin{scope}[scale=1.2] 
\fill (0,0) circle (3pt); 
\fill (1,0) circle (3pt); 
\fill (3,0) circle (3pt); 
\fill (4,0) circle (3pt); 
\fill (5,0) circle (3pt); 
\draw (0,0)--(1,0)--(1.5,0);
\draw [dashed] (1.5,0)--(2.5,0);
\draw (2.5,0)--(3,0)--(4,0);	
\draw (4,0.05)--(5,0.05);
\draw (4,-0.05)--(5,-0.05);
\draw (4.6,0.15)--(4.45,0);
\draw (4.45,0)--(4.6,-0.15);
\draw (0,-0.3) node {$s_1$};	
\draw (1,-0.3) node {$s_2$};	
\draw (3,-0.3) node {$s_{n-2}$};	
\draw (4,-0.3) node {$s_{n-1}$};	
\draw (5,-0.3) node {$s_n$};
\end{scope}
\end{tikzpicture}
&
$E_8$ 
&
\begin{tikzpicture}[scale=0.7,double distance=3pt,thick]
\begin{scope}[scale=1] 
\fill (0,-.5) circle (3pt); 
\fill (1,-.5) circle (3pt); 
\fill (2,-.5) circle (3pt); 
\fill (3,-.5) circle (3pt); 
\fill (4,-.5) circle (3pt); 
\fill (5,-.5) circle (3pt); 
\fill (6,-.5) circle (3pt); 
\fill (2,-1.5) circle (3pt); 
\draw (0,-.5)--(1,-.5)--(2,-.5)--(3,-.5)--(4,-.5)--(5,-.5)--(6,-.5);
\draw (2,-.5)--(2,-1.5);
\draw (0,-0.25) node {$s_1$};	
\draw (1,-0.25) node {$s_3$};
\draw (2,-0.25) node {$s_4$};
\draw (3,-0.25) node {$s_5$};
\draw (4,-0.25) node {$s_6$};
\draw (5,-0.25) node {$s_7$};
\draw (6,-0.25) node {$s_8$};
\draw (2.35,-1.4) node {$s_2$};	
\end{scope}
\end{tikzpicture}
 \\
\hline
$D_n,\ n \geq 4$ 
&
\begin{tikzpicture}[scale=0.7,double distance=3pt,thick]
\begin{scope}[scale=1.2] 
\fill (0,0) circle (3pt); 
\fill (1,0) circle (3pt); 
\fill (3,0) circle (3pt); 
\fill (4,0) circle (3pt); 
\fill (5,0.5) circle (3pt); 
\fill (5,-0.5) circle (3pt); 
\draw (0,0)--(1,0)--(1.5,0);
\draw [dashed] (1.5,0)--(2.5,0);
\draw (2.5,0)--(3,0)--(4,0); 
\draw (4,0)--(5,0.5);
\draw (4,0)--(5,-0.5);
\draw (0,-0.3) node {$s_1$};	
\draw (1,-0.3) node {$s_2$};	
\draw (3,-0.3) node {$s_{n-3}$};	
\draw (4,-0.3) node {$s_{n-2}$};	
\draw (5.2,0.8) node {$s_{n-1}$};
\draw (5,-0.8) node {$s_{n}$};
\end{scope}
\end{tikzpicture}
&
$F_4$ 
&
\begin{tikzpicture}[scale=0.7,double distance=3pt,thick]
\begin{scope}[scale=1.2] 
\fill (0,.5) circle (3pt); 
\fill (1,.5) circle (3pt); 
\fill (2,.5) circle (3pt); 
\fill (3,.5) circle (3pt); 
\draw (0,.5)--(1,.5);
\draw (2,.5)--(3,.5);
\draw (1,0.55)--(2,0.55);
\draw (1,0.45)--(2,0.45);
\draw (1.4,0.65)--(1.55,.5);
\draw (1.55,.5)--(1.4,0.35);
\draw (0,.15) node {$s_1$};	
\draw (1,.15) node {$s_2$};
\draw (2,.15) node {$s_3$};
\draw (3,.15) node {$s_4$};
\end{scope}
\end{tikzpicture}
 \\
\hline

&

&
$G_2$ 
&
\begin{tikzpicture}[scale=0.7,double distance=3pt,thick]
\begin{scope}[scale=1.2] 
\fill (0,-.25) circle (3pt); 
\fill (1,-.25) circle (3pt); 
\draw (0,-.18)--(1,-.18);
\draw (0,-.25)--(1,-.25);
\draw (0,-.32)--(1,-.32);
\draw (0.6,-.1)--(0.45,-.25);
\draw (0.45,-.25)--(0.6,-0.4);
\draw (0,-0.5) node {$s_1$};	
\draw (1,-0.5) node {$s_2$};
\end{scope}
\end{tikzpicture}
 \\
\hline
\end{tabular}
\caption{Dynkin diagrams for all types; labels consistent with \cite{Bour46} and \cite{sage}.}
\end{center}
\end{table}%

\subsection{Root systems}

The group $W$ acts by linear transformations on a real $n$-dimensional vector space $V$, which we may identify with $\R^n$. The vector space $V$ admits an ordered basis $\Delta = (\alpha_i)_{i \in [n]}$ of simple roots, and a symmetric bilinear form $B(\alpha_i, \alpha_j) = -\cos \frac{\pi}{m(i,j)}$, where $m(i,j)$ is the $(i,j)$-entry of the associated Coxeter matrix.  

Given any $s_i \in S$ and $\alpha_j \in \Delta$, the action of $W$ on $V$ is defined by 
\begin{equation}\label{eq:cartan}
s_i(\alpha_j) = \alpha_j - c_{ji} \alpha_i,
\end{equation}
where $c_{ji}$ is the $(j,i)$-entry of the associated Cartan matrix.  The set 
\[ \Phi = \{w(\alpha_i) \mid w \in W, \ i \in [n]\} \subset V\] 
forms the root system for $(W,S)$, and the elements are called roots. Given any roots $\alpha, \beta \in \Phi$, the action of $W$ on $V$ can also be expressed via
\begin{equation}\label{eq:rootreflect}
s_\alpha(\beta) = \beta - \frac{2(\alpha, \beta)}{(\alpha, \alpha)}\alpha,
\end{equation}
where $(\cdot, \cdot)$ denotes the Euclidean dot product. Let $\Phi^+$ (resp.~$\Phi^-$) denote the positive (resp.~negative) roots in $\Phi$, meaning those which have non-negative (resp.~non-positive) coefficients when expanded in terms of the basis $\Delta$. Let $\rho$ denote the half-sum of the positive roots in $\Phi^+$.

Given any $\alpha \in \Phi^+$, write $\alpha = w(\alpha_i)$ for some $w \in W$ and $\alpha_i \in \Delta$. The assignment 
\begin{equation}\label{eq:bij}
s_\alpha = ws_i w^{-1} 
\end{equation}
is a bijective correspondence between positive roots in $\Phi^+$ and reflections in $T$; see Proposition 4.5 in \cite{BB}. In particular,  we have $s_{\alpha_i} = s_i$. Moreover,  the element $ws_iw^{-1} \in T$ acts on $V$ as the reflection in the hyperplane orthogonal to the root vector $\alpha = w(\alpha_i) \in \Phi^+$.
For any $w\in W$, denote by
 \[N(w)=\{ \alpha\in \Phi^+\mid w(\alpha)\in \Phi^-\}.\]    
 Recall from \cite[Proposition 4.4.4]{BB} that $\ell(w)=|N(w)|$ for any $w \in W$. For any $\alpha \in \Phi^+$, we denote $N(\alpha)=N(s_{\alpha})$. In particular, we then have $\ell(s_\alpha) = |N(\alpha)|$. 
 
 Given any $i \in [n]$, define a linear functional $\alpha_i^\vee \in V^*$ by $\langle \alpha_i^\vee, v \rangle := 2B(\alpha_i, v)$ for any $v \in V$, where $\langle \cdot, \cdot \rangle: V^* \times V \to \Z$ denotes the evaluation pairing. The ordered set $\Delta^\vee = (\alpha_i^\vee)_{i \in [n]}$ is then a basis of simple coroots for the dual space $V^*$. 
For any root $\alpha = w(\alpha_i)$, the corresponding coroot $\alpha^\vee$ is defined by $\langle \alpha^\vee, v \rangle := 2B(\alpha, v)$ for $v \in V$.  Identifying $V$ and $V^*$, we may also write $\alpha^\vee = \frac{2\alpha}{(\alpha, \alpha)}$ for any root $\alpha \in \Phi$. The set  $\Phi^\vee = \{w(\alpha_i)^\vee  \mid w \in W, \ i \in [n]\} \subset V^*$ forms the system of coroots for $(W,S)$.  The height of the coroot $\alpha^\vee \in \Phi^\vee$ is the sum of the coefficients $c_i$ in the expansion $\alpha^\vee = \sum c_i \alpha_i^\vee$ via $\Delta^\vee$. Equivalently, the height of $\alpha^\vee$ equals $\height (\alpha^\vee) = \langle \alpha^\vee, \rho\rangle$.

\section{Classical Types}\label{sec:classical}

This section is dedicated to the proof of Theorem \ref{thm:main} in the classical types.
We refer the reader to the plates in \cite{Bour46}, which provide explicit formulas for the roots and Cartan matrices. Note that we correct a typo in item (II) in each of the plates for types $C_n$ and $D_n$ below.

\subsection{Type $A_n$}

 In type $A_n$, the roots are 
 \[ \Phi = \left\{ e_i-e_j\  \middle| \ 1\leq i,j\leq n+1 \right\}, \]
 where here $e_i$ denotes the $i^{\text{th}}$ standard basis vector in $\R^{n+1}$. We choose the simple roots to be those of the form $\Delta = \{\alpha_i\}= \{e_i-e_{i+1} \mid i \in [n] \}$, so that the positive roots are 
 \[\Phi^+=\{\alpha_{ij} := e_i-e_{j+1}\mid 1 \leq i<j \leq n\}.\]  
  In the geometric realization of $(W,S)$, the reflection $s_{\alpha_i}$ acts by interchanging coordinates $i$ and $i+1$ and is thus identified with the simple transposition $s_i \in S_{n+1} \cong W$. 
 
 Since $\alpha^\vee = \frac{2\alpha}{(\alpha, \alpha)}$ for any $\alpha \in \Phi$, we observe that $\alpha^\vee = \alpha$ for all $\alpha \in \Phi$ in type $A_n$.  The corresponding basis of simple coroots is thus $\Delta^\vee = \{\alpha_i^\vee = \alpha_i\}$, and we can write the positive coroot $\alpha_{ij}^\vee = (e_i-e_j)^\vee = \alpha_i^\vee + \cdots + \alpha_j^\vee$.  In particular, note that 
 $\height(\alpha_{ij}^\vee) = j-i+1$. Since $A_n$ is simply-laced, Lemma 4.3 of \cite{BFP} says that $\ell(s_\alpha) = \langle 2\rho, \alpha^\vee \rangle - 1 = 2\height(\alpha^\vee)-1$ for all $\alpha \in \Phi^+$. Therefore, $\ell(s_{\alpha_{ij}}) = 2(j-i+1)-1 = 2(j-i)+1$.
 
 For any fixed $\alpha_{ij} \in \Phi^+$, we claim that 
 \begin{equation}\label{eq:Areflect}
s_{e_i-e_{j+1}} =  s_{\alpha_{ij}} = s_is_{i+1} \cdots s_{j-1}s_j s_{j-1} \cdots s_{i+1}s_i 
 \end{equation}
 is a reduced expression for the reflection corresponding to $\alpha_{ij} = \alpha_i + \cdots + \alpha_j$. In type $A$, this fact is well-known, but we include a proof for the sake of completeness, since there are large type $A$ subsystems in the other types, and we thus rely inductively on this type $A$ result.  
 
 Note that the expression in \eqref{eq:Areflect} uses exactly $\ell(s_{\alpha_{ij}}) = 2(j-i)+1$ simple generators, and is thus reduced if it indeed represents $s_{\alpha_{ij}}$.  Let $w = s_is_{i+1} \cdots s_{j-1}$.  Using the bijective correspondence between positive roots and reflections from \eqref{eq:bij}, if we can prove that $w(\alpha_j) = \alpha_{ij}$, then $s_{\alpha_{ij}} = ws_j w^{-1}$ is indeed the reflection corresponding to $\alpha_{ij}$. By \eqref{eq:cartan}, in type $A_n$ we have
 \begin{equation}\label{eq:Aaction}
s_i (\alpha_j) = \begin{cases} - \alpha_j & \text{if $i=j$} \\ \alpha_i + \alpha_j & \text{if $|i -j| =1$} \\ 
\alpha_j & \text{otherwise}.
\end{cases}
\end{equation}
We can now directly compute via \eqref{eq:Aaction} that
\begin{equation*}
w(\alpha_j) = s_is_{i+1} \cdots s_{j-1}(\alpha_j) = \alpha_i + \cdots + \alpha_j = \alpha_{ij},
\end{equation*}
as required.  Therefore, the expressions in \eqref{eq:Areflect} form a system of palindromic reduced expressions for the $n \choose 2$ $= \frac{n(n-1)}{2}$ distinct reflections in $W$ of type $A_n$.

\subsection{Type $C_n$}\label{sec:Cn}

 In type $C_n$, the roots are 
 \[ \Phi = \left\{ \pm e_i\pm e_j\  \middle| \ 1\leq i < j \leq n \right\} \cup \{ \pm 2 e_i \mid i \in [n]\}, \]
 where here $e_i$ denotes the $i^{\text{th}}$ standard basis vector in $\R^{n}$. We choose the simple roots to be $\Delta = \{e_i-e_{i+1} \mid i \in [n-1] \} \cup \{ 2e_n\}$, where $\alpha_i = e_i-e_{i+1}$ as in type $A_{n-1}$, but $\alpha_n = 2e_n$ is distinct to type $C_n$.   
 The $n^2$ positive roots then come in three natural families:
 \[\Phi^+=\{e_i-e_{j}\mid 1 \leq i<j \leq n\} \cup \{ e_i+e_j \mid 1 \leq i < j \leq n \} \cup \{ 2e_i \mid i \in [n]\}.\]  
 We study each of these families separately.
 
 For the $n \choose 2$ $=\frac{n(n-1)}{2}$ positive roots in $\{e_i-e_{j}\mid 1 \leq i<j \leq n\}$, we can write 
 \[ e_i-e_j = \alpha_i + \alpha_{i+1} + \cdots + \alpha_{j-1},\]
which are the roots for the underlying type $A_{n-1}$ root system generated by the first $n-1$ nodes of  the Dynkin diagram.  Therefore, for any fixed $e_i-e_j \in \Phi^+$, we have the following palindromic reduced expression
 \begin{equation}\label{eq:C1reflect}
 s_{e_i-e_j} = s_is_{i+1} \cdots s_{j-2}s_{j-1} s_{j-2} \cdots s_{i+1}s_i 
 \end{equation}
by our proof in type $A$.
 
For any fixed $i \in [n]$, we claim that 
 \begin{equation}\label{eq:C2reflect}
 s_{2e_{i}} = s_is_{i+1} \cdots s_{n-1}s_n s_{n-1} \cdots s_{i+1}s_i 
 \end{equation}
 is a reduced expression for each of the reflections corresponding to the $n$ positive roots $2e_i$.  Let $w = s_i \cdots s_{n-1}$. We first show that $w(\alpha_n) = 2e_i$.   Write 
 \[ 2e_i = 2\alpha_i + \cdots + 2\alpha_{n-1} + \alpha_n,\]
  correcting a typo in item (II) of the type $C_n$ plate in \cite{Bour46}, which omits the 2s on the righthand side of this expression. By \eqref{eq:cartan}, in type $C_n$ we have 
  \begin{equation}\label{eq:Caction}
s_i (\alpha_j) = \begin{cases} - \alpha_j & \text{if $i=j$} \\ \alpha_i + \alpha_j & \text{if $|i -j| =1$ and $j \neq n$} \\ 
2\alpha_{i}+\alpha_j & \text{if $i=n-1$ and $j=n$} \\
\alpha_j & \text{otherwise}.
\end{cases}
\end{equation}
We can now directly compute via \eqref{eq:Caction} that
\begin{equation*}
w(\alpha_n) = s_is_{i+1} \cdots s_{n-1}(\alpha_n) = 2\alpha_i + \cdots + 2\alpha_{n-1} + \alpha_n = 2e_i,
\end{equation*}
as required to verify that $s_{2e_i} = ws_j w^{-1}$ via \eqref{eq:bij}.

Note that the expression in \eqref{eq:C2reflect} uses exactly $2(n-i)+1$ simple reflections. To see that this expression for $s_{2e_i}$ is reduced, we prove that $\ell(s_{2e_i}) = |N(2e_i)|= 2(n-i)+1$.  Recall that $i \in [n]$ is fixed here, and define the following subset of positive roots:
 \[\Phi^+_{2e_i}=\{e_i-e_{j}\mid i<j \leq n\} \cup \{ e_i+e_j \mid i < j \leq n \} \cup \{ 2e_i\}.\]  
 We claim that $\Phi^+_{2e_i} = N(2e_i)$. For each $\alpha \in \Phi^+_{2e_i}$, we use \eqref{eq:rootreflect} to directly calculate $s_{2e_i}(\alpha)= \alpha - (2e_i,\alpha)e_i$.  For any $\alpha = e_i-e_j$ with $i < j \leq n$, we have $(2e_i,e_i-e_j) = 2$ and so $s_{2e_i}(e_i-e_j)= (e_i-e_j) - 2e_i = -e_i-e_j \in \Phi^-$. Similarly, for any $\alpha = e_i+e_j$ with $i < j \leq n$, we have  $s_{2e_i}(e_i+e_j)= (e_i+e_j) - 2e_i = -e_i+e_j \in \Phi^-$.  Finally, $s_{2e_i}(2e_i) = -2e_i \in \Phi^-$, completing the proof that $\Phi^+_{2e_i} \subseteq N(2e_i)$. 
 
 Note that there are exactly $n-i$ elements in each of the first and second subsets of $\Phi^+_{2e_i}$, so that altogether $|\Phi^+_{2e_i}| =2(n-i)+1 \leq |N(2e_i)| = \ell(s_{2e_i})$. Since the expression  in \eqref{eq:C2reflect} uses exactly $2(n-i)+1$ simple reflections, then this inequality implies that the expression for the reflection $s_{2e_i}$ in \eqref{eq:C2reflect} is reduced for any $i \in [n]$. Consequently, we have $\ell(s_{2e_i}) = 2(n-i)+1$, in which case $\Phi^+_{2e_i} = N(2e_i)$.

 For the remaining $n \choose 2$ $=\frac{n(n-1)}{2}$ positive roots $\{e_i+e_{j}\mid 1 \leq i<j \leq n\}$, we claim that 
\begin{equation}\label{eq:C3reflect}
s_{e_i+e_j} = \left(s_js_{j+1} \cdots s_{n-1} \cdot s_is_{i+1} \cdots s_{n-2}\right)\left(s_n s_{n-1} s_n\right) \left( s_{n-2} \cdots s_{i+1}s_i \cdot s_{n-1} \cdots s_{j+1} s_j\right)
\end{equation}
is a reduced expression.  Here, the product $s_js_{j+1}\cdots s_{n-1}$ is trivial if $j=n$, and the product $s_i s_{i+1} \cdots s_{n-2}$ is trivial if $i=n-1$.

We first provide some context for the choice of conjugating element in \eqref{eq:C3reflect}, generalizing our discussion in Example \ref{ex:introC4}. Consider the subgroup of $W$ of type $A_{n-1}$ given by  $\langle s_1, \dots, s_{n-1}\rangle \cong S_n$, and the maximal parabolic subgroup $P$ obtained by omitting the reflection $s_{n-2}$.  The minimial length coset representatives in $S_n/P$ are  in bijection with the Young diagrams contained in an $(n-2)\times 2$ rectangle, and a corresponding reduced word can be read by overlaying the Young diagram in English notation on the $(n-2)\times 2$ grid labeled by:
\begin{equation*}\label{Young}
  \begin{tabular}{ | c | c |}
    \hline
    $n-2$ & $n-1$  \\ \hline
    $n-3$ & $n-2$  \\ \hline
    $n-4$ & $n-3$  \\ \hline
    $\vdots$ & $\vdots$  \\ \hline
    $3$ & $4$  \\ \hline
    $2$ & $3$  \\ \hline
    $1$ & $2$  \\ 
    \hline
  \end{tabular} \quad\quad\quad \rightsquigarrow \quad\quad \quad
  \tableau[Ys]{ \liz{7}& \tr{8} \\ \liz{6} & \tr{7} \\ \liz{5} & \tr{6} \\ \liz{4} \\ \liz{3}  \\ \liz{2} \\ \liz{1}\\ } \quad \longleftrightarrow \quad s_{\tr{6}}s_{\tr{7}}s_{\tr{8}}s_{\liz{1}}s_{\liz{2}}s_{\liz{3}}s_{\liz{4}}s_{\liz{5}}s_{\liz{6}}s_{\liz{7}}.
\end{equation*}
Our convention will be to build the word from left to right by reading the labels up each column of the diagram, moving from the right column to the left column, as shown in the $n=9$ example in the figure above. Denote by $w_{ij}$ the element corresponding to the Young diagram whose first column extends down to label $i$ and second column extends down to label $j$. If the second column is empty, define $j=n$, and if the first column is empty, define $i=n-1$ so that both products are trivial.  (In the example above, we have $i=1$ and $j=6$, and $w_{16} = s_6s_7s_8s_1s_2s_3s_4s_5s_6s_7$.) Explicitly, the conjugating element in \eqref{eq:C3reflect} is then
\[ w_{ij} = s_js_{j+1} \cdots s_{n-1} \cdot s_is_{i+1} \cdots s_{n-2}.\]

In this notation, we aim to prove that
\[
 s_{e_i+e_j} = w_{ij}(s_n s_{n-1} s_n) w_{ij}^{-1}
\]
is a reduced expression for the reflection corresponding to the positive root $e_i+e_j$.
Write 
  \[
e_i+e_j  =
\begin{cases}
\alpha_i + \cdots + \alpha_{j-1}+2\alpha_j + \cdots + 2\alpha_{n-1}+\alpha_n & \text{if $j \neq n$,}\\
\alpha_i + \cdots + \alpha_n & \text{if $j=n$}.
\end{cases}
\]
From \eqref{eq:Caction}, compute that $s_n(\alpha_{n-1}) = \alpha_{n-1}+\alpha_n = e_{n-1}+e_{n}$, and so $s_{e_{n-1}+e_{n}} = s_ns_{n-1}s_n$ by \eqref{eq:bij}, which is clearly reduced. We show that the roots $e_i+e_j$ correspond bijectively to the 
 conjugates of $s_{e_{n-1}+e_{n}}$ by the $n \choose n-2$ $=$ $n \choose 2$ distinct minimal length coset representatives of $S_n/P$. 
 
Consider the case where $j=n$.  Compute using \eqref{eq:Caction} that 
\[ w_{in}s_n(\alpha_{n-1}) = w_{in}(\alpha_{n-1}+\alpha_n) = s_is_{i+1} \cdots s_{n-2}(\alpha_{n-1}+\alpha_n) = \alpha_i + \cdots +\alpha_n = e_i+e_n, \]
verifying \eqref{eq:C3reflect} for $j=n$ via \eqref{eq:bij}.
Now consider the case where the second column of the Young diagram for $w_{ij}$ is nonempty, equivalently $j \leq n-1$.  Compute using \eqref{eq:Caction} and the $j=n$ case that
\[ w_{ij}s_n(\alpha_{n-1}) = w_{ij}(\alpha_{n-1}+\alpha_n) =  s_js_{j+1} \cdots s_{n-1} \cdot s_is_{i+1} \cdots s_{n-2}(\alpha_{n-1}+\alpha_n)= \phantom{XXXX}\]
\[  =s_js_{j+1} \cdots s_{n-1} ( \alpha_i + \cdots +\alpha_n)  = \alpha_i + \cdots + \alpha_{j-1}+2\alpha_j + \cdots + 2\alpha_{n-1}+\alpha_n = e_i+e_j, \]
verifying \eqref{eq:C3reflect} for $j \neq n$ via \eqref{eq:bij}.  We have proved that $s_{e_i+e_j} = w_{ij}(s_n s_{n-1} s_n) w_{ij}^{-1}$ for all $1 \leq i < j \leq n$.  

Note that the expression in \eqref{eq:C3reflect} uses exactly $4n - 2(i+j)+1$ simple reflections. To see that this expression for $s_{e_i+e_j}$ is reduced, we prove that $\ell(s_{e_i+e_j}) = |N(e_i+e_j)| = 4n - 2(i+j)+1$.
 Define the following subset of positive roots:
 \begin{align*}
 \Phi^+_{e_i+e_j} &=\{e_k-e_l\mid i=k<l\neq j, \text{ or } j=k<l\} \\ 
 & \quad \cup  \{e_k+e_l\mid i=k \text{ and } j\leq l, \text{ or }j=l \text{ and } i<k \neq j \} \\
 & \quad \cup\{2e_i,2e_j\}.
 \end{align*}
 We claim that $\Phi^+_{e_i+e_j} = N(e_i+e_j)$.  For each $\alpha \in \Phi^+_{e_i+e_j}$, we use \eqref{eq:rootreflect} to directly calculate $s_{e_i+e_j}(\alpha)= \alpha - (e_i+e_j,\alpha)(e_i+e_j)$.  For any $\alpha =e_k-e_l$ such that $i = k< l \neq j$, we have $(e_i+e_j, e_i-e_l) = 1$ and so $s_{e_i+e_j}(e_i-e_l) = (e_i-e_l) - (e_i+e_j) =  -e_l-e_j \in \Phi^-$.  Note that there are exactly $n-i-1$ positive roots of the form $e_i-e_l$ with $l \neq j$.  Similarly, for any $\alpha = e_k-e_l$ such that $j=k<l$, we have $(e_i+e_j, e_j-e_l) = 1$ and so $s_{e_i+e_j}(e_j-e_l) = (e_j-e_l) - (e_i+e_j) =  -e_i-e_l \in \Phi^-$.  Note that there are exactly $n-j$ positive roots of the form $e_j-e_l$.  Altogether, we have identified $(n-i-1)+(n-j) = 2n-i-j-1$ positive roots of the form $e_k-e_l$ in $N(e_i+e_j)$.
 
 Now consider $\alpha = e_k+e_l$ such that $i=k$ and $j \leq l$.  If $\alpha = e_i+e_j$, then clearly $s_{e_i+e_j}(e_i+e_j) = -e_i-e_j \in \Phi^-$.  If $j <l$, then $(e_i+e_j,e_i+e_l) = 1$ and so $s_{e_i+e_j}(e_i+e_l) = -e_j+e_l \in \Phi^-$ since $j<l$. Note that there are exactly $n-j+1$ positive roots of the form $e_i+e_l$ such that $j \leq l$.  Now consider $\alpha = e_k+e_l$ such that $j=l$ and $i<k \neq j$. Then $(e_i+e_j, e_k+e_j) = 1$ and so $s_{e_i+e_j}(e_k+e_j) = -e_i+e_k \in \Phi^-$ since $i < k$. Note that there are exactly $n-i-1$ positive roots of the form $e_k+e_j$ such that $i<k\neq j$. Altogether, we have identified $(n-j+1)+(n-i-1) =2n-i-j $ positive roots of the form $e_k+e_l$ in $N(e_i+e_j)$. 
 
 Finally, $s_{e_i+e_j}(2e_i) = -2e_j \in \Phi^-$ and $s_{e_i+e_j}(2e_j) = -2e_i \in \Phi^-$, which concludes our proof that $\Phi^+_{e_i+e_j} \subseteq N(e_i+e_j)$.  We thus have $|\Phi^+_{e_i+e_j}| = (2n-i-j-1)+(2n-i-j)+2 = 4n-2(i+j)+1 \leq |N(e_i+e_j)| = \ell(s_{e_i+e_j})$. Since the expression  in \eqref{eq:C3reflect} uses exactly $4n-2(i+j)+1$ simple reflections, then this inequality implies that the expression for the reflection $s_{e_i+e_j}$ in \eqref{eq:C3reflect} is reduced for any $1 \leq i < j \leq n$. Consequently, we have $\ell(s_{e_i+e_j}) = 4n-2(i+j)+1$, in which case $\Phi^+_{e_i+e_j} = N(e_i+e_j)$.

\subsection{Type $B_n$}

The Weyl groups in types $B_n$ and $C_n$ are identical, and therefore so is the set of reflections $T = \{wsw^{-1} \mid w \in W, s \in S\}$.  However, the roles of the roots and coroots are reversed in types $B_n$ and $C_n$, and so the explicit correspondence between reduced words for reflections and positive roots requires minor relabeling in type $B_n$.

 In type $B_n$, the roots are 
 \[ \Phi = \left\{ \pm e_i\pm e_j\  \middle| \ 1\leq i < j \leq n \right\} \cup \{ \pm  e_i \mid i \in [n]\}, \]
 where $e_i$ denotes the $i^{\text{th}}$ standard basis vector in $\R^{n}$. We choose the simple roots to be $\Delta = \{e_i-e_{i+1} \mid i \in [n-1] \} \cup \{ e_n\}$, where $\alpha_i = e_i-e_{i+1}$ as in type $A_{n-1}$, but here $\alpha_n = e_n$ is distinct to type $B_n$.  The $n^2$ positive roots are then
 \[ \Phi^+ = \{ e_i\pm e_j \mid 1 \leq i < j \leq n\} \cup \{ e_i \mid i \in [n]\}. \]
Compute directly that $(e_i \pm e_j)^\vee = e_i \pm e_j$, whereas $e_i^\vee = 2e_i$.  Therefore, the reflections corresponding to the roots $e_i \pm e_j \in \Phi^+$ are labeled identically to those in type $C_n$, whereas the reflections $s_{2e_i}$ in type $C_n$ are relabeled by $s_{e_i}$ in type $B_n$. We thus have the following palindromic reduced expressions
\begin{align*}
s_{e_i-e_j} & = s_is_{i+1} \cdots s_{j-2}s_{j-1} s_{j-2} \cdots s_{i+1}s_i \nonumber \\
s_{e_i} & = s_is_{i+1} \cdots s_{n-1}s_n s_{n-1} \cdots s_{i+1}s_i \\
s_{e_i+e_j} & = \left(s_js_{j+1} \cdots s_{n-1} \cdot s_is_{i+1} \cdots s_{n-2}\right)\left(s_n s_{n-1} s_n\right) \left( s_{n-2} \cdots s_{i+1}s_i \cdot s_{n-1} \cdots s_{j+1} s_j\right) \nonumber
\end{align*}
for the respective $n \choose 2$ $+ \; n\; +$ $n \choose 2$ $=n^2$ reflections in type $B_n$.

\subsection{Type $D_n$}

 In type $D_n$, the roots are 
 \[ \Phi = \left\{ \pm e_i\pm e_j\  \middle| \ 1\leq i < j \leq n \right\}, \]
 where $e_i$ denotes the $i^{\text{th}}$ standard basis vector in $\R^{n}$.
  We choose the simple roots to be $\Delta =  \{e_i-e_{i+1} \mid i \in [n-1] \} \cup \{ e_{n-1}+e_n\}$, where $\alpha_i = e_i-e_{i+1}$ as in all other classical types, but here $\alpha_n = e_{n-1}+e_n$ is distinct to type $D_n$.    The $n^2-n$ positive roots are then
 \[\Phi^+=\{e_i-e_{j}\mid 1 \leq i<j \leq n\} \cup \{ e_i+e_j \mid 1 \leq i < j \leq n \}.\]

 For the $n \choose 2$ $=\frac{n(n-1)}{2}$ positive roots in $\{e_i-e_{j}\mid 1 \leq i<j \leq n\}$, we can write 
 \[ e_i-e_j = \alpha_i + \alpha_{i+1} + \cdots + \alpha_{j-1}\]
which are the roots for the underlying type $A_{n-1}$ root system generated by the first $n-1$ nodes of the Dynkin diagram.  (Note that this corrects a typo in item (II) of the type $D_n$ plate in \cite{Bour46}, which omits $\alpha_i$ from this expansion.) Therefore, for any fixed $e_i-e_j \in \Phi^+$, we have the following palindromic reduced expression
 \begin{equation}\label{eq:D1reflect}
 s_{e_i-e_j} = s_is_{i+1} \cdots s_{j-2}s_{j-1} s_{j-2} \cdots s_{i+1}s_i 
 \end{equation}
of length $2(j-i)-1$, by our proof in type $A$.

Next we consider the type $A_{n-1}$ subsystem generated by $\{s_1, \dots, s_{n-2}, s_n\}$.  The reflections which do not use $s_n$ in any reduced expression already appear among the reflections $s_{e_i-e_j}$ expanded in \eqref{eq:D1reflect}. For the remaining $n-1$ reflections in this system, 
we claim that 
 \begin{equation}\label{eq:D2reflect}
 s_{e_i+e_n} = 
 s_is_{i+1} \cdots s_{n-2}s_n s_{n-2} \cdots s_{i+1}s_i
 \end{equation}
 is a reduced expression for all $1 \leq i < n$, where the product $s_is_{i+1}\cdots s_{n-2}$ is trivial if $i = n-1$.  In case $i=n-1$, recall that $\alpha_n = e_{n-1}+e_n$ is a simple reflection, and so $s_{\alpha_n} = s_n$ by definition.  For $1 \leq i \leq n-2$, compute using \eqref{eq:cartan} that in type $D_n$
 \[ s_i s_{i+1} \cdots s_{n-2}(\alpha_n) = \alpha_i + \alpha_{i+1} + \cdots + \alpha_{n-2} + \alpha_n = e_i+e_n.\]
 Therefore, \eqref{eq:D2reflect} provides expressions for each of the $n-1$ reflections $s_{e_i+e_n}$ by \eqref{eq:bij}.
 
 To prove that the expressions in \eqref{eq:D2reflect} are reduced, we observe that $\alpha^\vee = \alpha$ for all $\alpha \in \Phi$ since type $D_n$ is simply-laced.  The basis of simple coroots is thus $\Delta^\vee = \{\alpha_i^\vee = \alpha_i\}$.  In particular, \[(e_i+e_n)^\vee = (\alpha_i + \alpha_{i+1} + \cdots + \alpha_{n-2} + \alpha_n)^\vee = \alpha_i^\vee + \alpha_{i+1}^\vee + \cdots + \alpha_{n-2}^\vee+ \alpha_n^\vee,\]
and so $\height(e_i+e_n)^\vee = n-i$. Since $D_n$ is simply-laced, Lemma 4.3 of \cite{BFP} says that $\ell(s_\alpha)= 2\height(\alpha^\vee)-1$ for all $\alpha \in \Phi^+$. Therefore, $\ell(s_{e_i+e_n}) = 2 \height(e_i+e_n)^\vee -1 = 2(n-i)-1$ for all $1 \leq i \leq n-1$.  Since the expression for $s_{e_i+e_n}$ in \eqref{eq:D2reflect} uses exactly $2(n-i)-1$ simple reflections, it must be reduced.
 
  For the remaining $n-1 \choose 2$ $=\frac{(n-1)(n-2)}{2}$ positive roots $\{e_i+e_{j}\mid 1 \leq i<j <n\}$, we claim that 
\begin{equation}\label{eq:D3reflect}
s_{e_i+e_j} = \left(s_js_{j+1} \cdots s_{n-2} \cdot s_i s_{i+1}\cdots s_{n-3}\right)\left(s_{n-1}s_{n-2}s_ns_{n-2}s_{n-1}\right) \left( s_{n-3} \cdots s_i \cdot s_{n-2} \cdots s_j\right)
\end{equation}
 is a reduced expression.  Here, the product $s_js_{j+1}\cdots s_{n-2}$ is trivial if $j=n-1$, and the product $s_i s_{i+1} \cdots s_{n-3}$ is trivial if $i=n-2$. To contextualize this conjugating element, consider the subgroup of $W$ of type $A_{n-2}$ generated by  $\langle s_1, \dots, s_{n-2} \rangle \cong S_{n-1}$, and the maximal parabolic subgroup $P$ obtained by omitting  $s_{n-3}$.  The minimial length coset representatives in $S_{n-1}/P$ are  in bijection with the Young diagrams contained in the $(n-3)\times 2$ rectangle labeled by
\begin{equation*}\label{Young}
  \begin{tabular}{ | c | c |}
    \hline
    $n-3$ & $n-2$  \\ \hline
    $n-4$ & $n-3$  \\ \hline
    $\vdots$ & $\vdots$  \\ \hline
    $2$ & $3$  \\ \hline
    $1$ & $2$  \\ 
    \hline
  \end{tabular} 
\end{equation*}
and read in the same manner as in type $C_{n-1}$. We show that the remaining roots $e_i+e_j$ in type $D_n$ correspond to the  conjugates of the reflection $s_{n-1}s_{n-2}s_ns_{n-2}s_{n-1}$ by the elements of the quotient $S_{n-1}/P$.  

Denote the conjugating element in \eqref{eq:D3reflect} by 
\[ w_{ij} = s_js_{j+1} \cdots s_{n-2} \cdot s_is_{i+1} \cdots s_{n-3}.\]
We aim to prove that
$ s_{e_i+e_j} = w_{ij}(s_{n-1}s_{n-2}s_ns_{n-2}s_{n-1}) w_{ij}^{-1}$
is a reduced expression for the reflection corresponding to the positive root $e_i+e_j$.
First consider the case $i = n-2$ and $j=n-1$, equivalently $w_{ij} =1$. Since the root subsystem spanned by $\alpha_{n-1}, \alpha_{n-2}, \alpha_n$ is type $A_3$, with an ordering on the roots which reverses the roles of $\alpha_{n-2}$ and $\alpha_{n-1}$, we know from our type $A$ results that the root $e_{n-2}+e_{n-1} = \alpha_{n-2}+\alpha_{n-1}+\alpha_n$ corresponds to the reflection $s_{n-1}s_{n-2}s_ns_{n-2}s_{n-1}$.  For $1 \leq i < j < n$, we have 
  \[
e_i+e_j  =
\begin{cases}
\alpha_i + \cdots + \alpha_{j-1}+2\alpha_j + \cdots + 2\alpha_{n-2}+\alpha_{n-1}+\alpha_n & \text{if $j < n-1$,}\\
\alpha_i + \cdots + \alpha_n & \text{if $j=n-1$}.
\end{cases}
\]
From the entries of the Cartan matrix in type $D_n$, we have from \eqref{eq:cartan} that
 \begin{equation}\label{eq:Daction}
 s_{n-2}(\alpha_n) = \alpha_{n-2}+\alpha_n, \quad
  s_{n-1}(\alpha_n) = \alpha_n, \quad
   s_{n}(\alpha_{n-2}) = \alpha_{n-2}+\alpha_n, \quad
    s_{n}(\alpha_{n-1}) = \alpha_{n-1},
 \end{equation}
and otherwise, the action $s_i(\alpha_j)$ agrees with the type $A_n$ formula from \eqref{eq:Aaction}. First consider the case $j=n-1$.  Compute using \eqref{eq:Aaction} and \eqref{eq:Daction} that 
\[ w_{i,n-1}(\alpha_{n-2}+\alpha_{n-1}+\alpha_n) = s_is_{i+1} \cdots s_{n-3}(\alpha_{n-2}+\alpha_{n-1}+\alpha_n) = \alpha_i + \cdots +\alpha_n = e_i+e_{n-1}, \]
verifying \eqref{eq:D3reflect} for $j=n-1$ by \eqref{eq:bij}. Finally, consider the case where the second column of the Young diagram for $w_{ij}$ is nonempty, equivalently $j <n-1$.  Compute using \eqref{eq:Aaction} and \eqref{eq:Daction}, as well as the $j=n-1$ case, that
\[ w_{ij}(\alpha_{n-2}+\alpha_{n-1}+\alpha_n) =  s_js_{j+1} \cdots s_{n-2} \cdot s_is_{i+1} \cdots s_{n-3}(\alpha_{n-2}+\alpha_{n-1}+\alpha_n) = \phantom{XXXXX}\]
\[  = s_js_{j+1} \cdots s_{n-2} ( \alpha_i + \cdots +\alpha_n)  =\alpha_i + \cdots + \alpha_{j-1}+2\alpha_j + \cdots + 2\alpha_{n-2}+\alpha_{n-1}+\alpha_n = e_i+e_j, \]
 verifying \eqref{eq:D3reflect} for $j < n-1$ by \eqref{eq:bij}.
  We have thus proven that $s_{e_i+e_j} = w_{ij}(s_{n-1}s_{n-2}s_ns_{n-2}s_{n-1}) w_{ij}^{-1}$ for all $1 \leq i < j < n$.  

To prove that the expressions in \eqref{eq:D3reflect} are reduced, we again use the fact that type $D_n$ is simply-laced.   Since $\alpha^\vee = \alpha$ for all $\alpha \in \Phi$, we have
  \[
(e_i+e_j)^\vee  =
\begin{cases}
\alpha_i^\vee + \cdots + \alpha_{j-1}^\vee+2\alpha_j^\vee + \cdots + 2\alpha_{n-2}^\vee+\alpha_{n-1}^\vee+\alpha_n^\vee & \text{if $j < n-1$,}\\
\alpha_i^\vee + \cdots + \alpha_n^\vee & \text{if $j=n-1$}.
\end{cases}
\]
and in particular, 
$ \height (e_i+e_j)^\vee  = 2n-i-j$ for all $1 \leq i < j < n$. By Lemma 4.3 of \cite{BFP}, we have $\ell(s_{e_i+e_j}) = 2 \height(e_i+e_j)^\vee -1 = 4n-2(i+j)-1 $ for all $1 \leq i < j < n$.  Since the expression for $s_{e_i+e_j}$ in \eqref{eq:D3reflect} uses exactly $2[(n-1-j)+(n-2-i)]+5 = 4n-2(i+j)-1$ simple reflections for all $1 \leq i < j < n$, it must be reduced. 

This concludes our proof of Theorem \ref{thm:main} in the classical types. \qed

\section{Exceptional Types}\label{sec:exceptional}

This section is dedicated to the proof of Theorem \ref{thm:main} in the exceptional types.
We discuss the exceptional types in order of increasing rank. We refer the reader to the plates in \cite{Bour46}, which provide explicit formulas for the roots and Cartan matrices. Note that we correct a typo in item (II) for type $G_2$, as well as two of the positive root expansions in item (II) of type $E_8$.

Beyond type $G_2$, we rely heavily upon our results in the classical types, as each exceptional group contains several parabolic subgroups of type $A, B, C$ and/or $D$. Throughout our discussion of the exceptional types, we typically denote $s_is_js_k = s_{ijk}$ for brevity.

\begin{table}[h]
\centering
  \begin{tabular}{ | c | c | c |}
    \hline
   Root $\alpha \in \Phi^+$ & Expand $\alpha$ in $\Delta$ & Reflection $s_\alpha$ \\ \hline\hline
   $\varepsilon_1-\varepsilon_2$  & $\alpha_1$ & $s_1$  \\ \hline
    $-2\varepsilon_1+\varepsilon_2+\varepsilon_3$  & $\alpha_2$ & $s_2$  \\ \hline
      $-\varepsilon_1+\varepsilon_3$  & $\alpha_1+\alpha_2$ & $s_{212} $ \\ \hline
        $\varepsilon_1-2\varepsilon_2+\varepsilon_3$  & $3\alpha_1+\alpha_2$ & $s_{121} $ \\ \hline
            $-\varepsilon_2+\varepsilon_3$  & $2\alpha_1+\alpha_2$ & $s_{12121} $ \\ \hline
               $-\varepsilon_1-\varepsilon_2+2\varepsilon_3$   & $3\alpha_1+2\alpha_2$ & $s_{21212} $ \\ \hline
      \end{tabular} 
  \vskip 5pt
\caption{Palindromic reduced expressions for the 6 reflections in type $G_2$.}\label{tab:G2}  
  \end{table}

\subsection{Type $G_2$}

In type $G_2$, we choose the basis of simple roots $\alpha_1 = \varepsilon_1 - \varepsilon_2, \alpha_2 = -2\varepsilon_1+\varepsilon_2+\varepsilon_3$, where $\varepsilon_i$ denotes the projection of the standard basis vector $e_i \in \R^3$ onto the hyperplane $\{ (x_1, x_2, x_3) \mid x_1+x_2+x_3 = 0\}$.  The positive roots are then
\[ \Phi^+ = \{ \alpha_1, \alpha_2, \alpha_1+\alpha_2, 2\alpha_1+\alpha_2, 3\alpha_1+\alpha_2, 3\alpha_1+2\alpha_2\}, \]
 where we point out that item (II) in the $G_2$ plate of \cite{Bour46} omits $\alpha_2$ from $\Phi^+$.  Table \ref{tab:G2} provides the expansion for each positive root in terms of $\{ \varepsilon_1, \varepsilon_2, \varepsilon_3\}$.

 Using the Cartan matrix in type $G_2$ together with \eqref{eq:cartan}, we have
 \[ s_1(\alpha_2) =3\alpha_1+\alpha_2 \quad \phantom{X} \quad s_2(\alpha_1) = \alpha_1+\alpha_2\]
  \[ s_1s_2(\alpha_1) =2\alpha_1+\alpha_2 \quad \quad s_2s_1(\alpha_2) = 3\alpha_1+2\alpha_2.\]
  We thus immediately obtain the palindromic expressions in Table \ref{tab:G2} via \eqref{eq:bij}, all of which are reduced by inspection.

\begin{table}[h]
\centering
  \begin{tabular}{ | c | c | c |}
    \hline
   Root $\alpha \in \Phi^+$ & Expand $\alpha$ in $\Delta$ & Reflection $s_\alpha$ \\ \hline\hline
    $e_2-e_3$ & $\alpha_1$ & $s_1$  \\ \hline
    $e_3-e_4$ & $\alpha_2$ & $s_2$  \\ \hline
    $e_2-e_4$ & $\alpha_1+\alpha_2$ & $s_{121}$  \\ \hline
    $e_2$ & $\alpha_1+\alpha_2+\alpha_3$ & $s_{12321}$  \\ \hline
    $e_3$ & $\alpha_2+\alpha_3$ & $s_{232}$  \\ \hline
    $e_4$ & $\alpha_3$ & $s_3$  \\ \hline

        $e_3+e_4$ & $\alpha_2+2\alpha_3$ & $s_{323}$  \\ \hline

    $e_2+e_4$ & $\alpha_1+\alpha_2+2\alpha_3$ & $s_{13231}$  \\ \hline 
            $e_2+e_3$ & $\alpha_1+2\alpha_2+2\alpha_3$ & $s_{2132312}$  \\ \hline\hline
    $\frac{1}{2}(e_1- e_2- e_3 - e_4 )$ & $\alpha_4$ & $s_4$  \\ \hline
     $\frac{1}{2}(e_1- e_2- e_3 +e_4 )$ & $\alpha_3+\alpha_4$ & $s_{434}$  \\ \hline
     $\frac{1}{2}(e_1- e_2+e_3 - e_4 )$ & $\alpha_2+\alpha_3+\alpha_4$ & $s_{42324}$  \\ \hline

      $\frac{1}{2}(e_1- e_2+ e_3 +e_4 )$ & $\alpha_2+2\alpha_3+\alpha_4$ & $s_{3423243}$  \\ \hline
           $e_1-e_2$ & $\alpha_2+2\alpha_3+2\alpha_4$ & $s_{43234}$  \\ \hline\hline
      $\frac{1}{2}(e_1+ e_2- e_3 - e_4 )$ & $\alpha_1+\alpha_2+\alpha_3+\alpha_4$ & $s_{1423241}$  \\ \hline
            $\frac{1}{2}(e_1+e_2- e_3 + e_4 )$ & $\alpha_1+\alpha_2+2\alpha_3+\alpha_4$ & $s_{134232431}$  \\ \hline  
                  $\frac{1}{2}(e_1+e_2+ e_3 - e_4 )$ & $\alpha_1+2\alpha_2+2\alpha_3+\alpha_4$ & $s_{21342324312}$  \\ \hline
                       $\frac{1}{2}(e_1+e_2+ e_3 + e_4 )$ & $\alpha_1+2\alpha_2+3\alpha_3+\alpha_4$ & $s_{3213423243123}$  \\ \hline

     $e_1-e_3$ & $\alpha_1+\alpha_2+2\alpha_3+2\alpha_4$ & $s_{1432341}$  \\ \hline
  
      $e_1-e_4$ & $\alpha_1+2\alpha_2+2\alpha_3+2\alpha_4$ & $s_{214323412}$  \\ \hline

    $e_1+e_4$ & $\alpha_1+2\alpha_2+4\alpha_3+2\alpha_4$ & $s_{32143234123}$  \\ \hline

   $e_1+e_3$  & $\alpha_1+3\alpha_2+4\alpha_3+2\alpha_4$ & $s_{2321432341232}$  \\ \hline
     $e_1$ & $\alpha_1+2\alpha_2+3\alpha_3+2\alpha_4$ & $s_{432134232431234}$  \\ \hline
    $e_1+e_2$ & $2\alpha_1+3\alpha_2+4\alpha_3+2\alpha_4$ & $s_{123214323412321}$  \\ \hline
      \end{tabular} 
  \vskip 5pt
\caption{Palindromic reduced expressions for the 24 reflections in type $F_4$.}\label{tab:F4}  
  \end{table}

\subsection{Type $F_4$}

In type $F_4$, we choose the simple roots $\Delta = \{\alpha_i\}$ to be
$\alpha_1 = e_2-e_3, \alpha_2 = e_3-e_4,  \alpha_3 = e_4,  \alpha_4 = \frac{1}{2}(e_1-e_2-e_3-e_4).$
The positive roots are then 
\[ \Phi^+ = \{ e_i \mid 1 \leq i \leq 4 \} \cup \{ e_i\pm e_j \mid 1 \leq i < j \leq 4 \} \cup \left\{ \frac{1}{2}(e_1\pm e_2\pm e_3 \pm e_4 )\right\}.\]
Table \ref{tab:F4} provides the expansion for each positive root in terms of $\Delta$.

It is clear from the Dynkin diagram and the choice of $\alpha_1, \alpha_2, \alpha_3$ that $W/ \langle s_4 \rangle$ has type $B_3$.  We thus obtain the first 9 entries in Table \ref{tab:F4}, by simply reindexing the standard basis vectors $e_i \mapsto e_{i+1}$ in type $B_3$.  In addition, $W/ \langle s_1 \rangle$ has type $C_3$, though the first and last labels are swapped compared to our labeling in  type $C$. Noting that the 4 reflections $s_2, s_3, s_{232}, s_{323}$ are common to $W / \langle s_1 \rangle$ and $W / \langle s_4 \rangle$, we obtain 5 distinct type $C_3$ reflections, which appear next in Table \ref{tab:F4}.  The remaining 10 reflections are distinct to type $F_4$, and are obtained by successively conjugating the 3 longest type $C_3$ reflections by simple reflections which increase the length. 

\begin{table}[h]
\centering
  \begin{tabular}{ | c | c | c |}
    \hline
  Root $\alpha \in \Phi^+$  & Expand $\alpha$ in $\Delta$ & Reflection $s_\alpha$ \\ \hline\hline
    01111 & $\alpha_1$ & $s_1$  \\ \hline
  10111   & $\alpha_1+\alpha_3$ & $s_{131}$  \\ \hline
    11011 & $\alpha_1+\alpha_3+\alpha_4$ & $s_{13431}$  \\ \hline
    00011  & $\alpha_1+\alpha_2+\alpha_3+\alpha_4$ & $s_{1342431}$  \\ \hline
    $-\varepsilon_1+\varepsilon_2$ & $\alpha_3$ & $s_3$  \\ \hline
    $-\varepsilon_1+\varepsilon_3$& $\alpha_3+\alpha_4$ & $s_{343}$  \\ \hline

     $-\varepsilon_2+\varepsilon_3$& $\alpha_4$ & $s_4$  \\ \hline
         $\varepsilon_1+\varepsilon_2$ & $\alpha_2$ & $s_2$  \\ \hline
              $\varepsilon_1+\varepsilon_3$ & $\alpha_2+\alpha_4$ & $s_{424}$  \\ \hline
          $\varepsilon_2+\varepsilon_3$ & $\alpha_2+\alpha_3+\alpha_4$ & $s_{34243}$  \\ \hline
     $-\varepsilon_3+\varepsilon_4$ & $\alpha_5$ & $s_5$  \\ \hline
         $-\varepsilon_2+\varepsilon_4$ & $\alpha_4+\alpha_5$ & $s_{454}$  \\ \hline
   $-\varepsilon_1+\varepsilon_4$ & $\alpha_3+\alpha_4+\alpha_5$ & $s_{34543}$  \\ \hline
    11101 & $\alpha_1+\alpha_3+\alpha_4+\alpha_5$ & $s_{1345431}$  \\ \hline

   00101  & $\alpha_1+\alpha_2+\alpha_3+\alpha_4+\alpha_5$ & $s_{132454231}$  \\ \hline
      01001  & $\alpha_1+\alpha_2+\alpha_3+2\alpha_4+\alpha_5$ & $s_{41324542314}$  \\ \hline
   10001  & $\alpha_1+\alpha_2+2\alpha_3+2\alpha_4+\alpha_5$ & $s_{3413245423143}$  \\ \hline

        $\varepsilon_1+\varepsilon_4$ & $\alpha_2+\alpha_4+\alpha_5$ & $s_{24542}$  \\ \hline
      $\varepsilon_2+\varepsilon_4$ & $\alpha_2+\alpha_3+\alpha_4+\alpha_5$ & $s_{3245423}$  \\ \hline
    $\varepsilon_3+\varepsilon_4$ & $\alpha_2+\alpha_3+2\alpha_4+\alpha_5$ & $s_{432454234}$  \\ \hline\hline
      $-\varepsilon_4+\varepsilon_5$  & $\alpha_6$ & $s_6$  \\ \hline 
        $-\varepsilon_3+\varepsilon_5$ & $\alpha_5+\alpha_6$ & $s_{656}$  \\ \hline 
        $-\varepsilon_2+\varepsilon_5$ & $\alpha_4+\alpha_5+\alpha_6$ & $s_{65456}$  \\ \hline

        $-\varepsilon_1+\varepsilon_5$ & $\alpha_3+\alpha_4+\alpha_5+\alpha_6$ & $s_{6543456}$  \\ \hline
           $\varepsilon_1+\varepsilon_5$& $\alpha_2+\alpha_4+\alpha_5+\alpha_6$ & $s_{6542456}$  \\ \hline
    $\varepsilon_2+\varepsilon_5$& $\alpha_2+\alpha_3+\alpha_4+\alpha_5+\alpha_6$ & $s_{652434256}$  \\ \hline
   $\varepsilon_3+\varepsilon_5$  & $\alpha_2+\alpha_3+2\alpha_4+\alpha_5+\alpha_6$ & $s_{46524342564}$  \\ \hline
     $\varepsilon_4+\varepsilon_5$  & $\alpha_2+\alpha_3+2\alpha_4+2\alpha_5+\alpha_6$ & $s_{5465243425645}$  \\ \hline  \hline
    11110  & $\alpha_1+\alpha_3+\alpha_4+\alpha_5+\alpha_6$ & $s_{134565431}$  \\ \hline\hline
      00110   & $\alpha_1+\alpha_2+\alpha_3+\alpha_4+\alpha_5+\alpha_6$ & $s_{16524342561}$ \\ \hline
        01010  & $\alpha_1+\alpha_2+\alpha_3+2\alpha_4+\alpha_5+\alpha_6$ & $s_{1465243425641}$  \\ \hline
          01100  & $\alpha_1+\alpha_2+\alpha_3+2\alpha_4+2\alpha_5+\alpha_6$ & $s_{154652434256451}$  \\ \hline
         10010    & $\alpha_1+\alpha_2+2\alpha_3+2\alpha_4+\alpha_5+\alpha_6$ & $s_{314652434256413}$  \\ \hline
                   10100      & $\alpha_1+\alpha_2+2\alpha_3+2\alpha_4+2\alpha_5+\alpha_6$ & $s_{31546524342564513}$  \\ \hline
                      11000   & $\alpha_1+\alpha_2+2\alpha_3+3\alpha_4+2\alpha_5+\alpha_6$ & $s_{4315465243425645134}$  \\ \hline
                     00000    & $\alpha_1+2\alpha_2+2\alpha_3+3\alpha_4+2\alpha_5+\alpha_6$ & $s_{243154652434256451342}$  \\ \hline
      \end{tabular} 
  \vskip 5pt
\caption{Palindromic reduced expressions for the 36 reflections in type $E_6$.}\label{tab:E6}  
  \end{table}
  
  \subsection{Type $E_6$}

In type $E_6$, let $\varepsilon_i$ denote the projection of the standard basis vector $e_i \in \R^8$ onto the subspace $\{ (x_1, \dots, x_8) \mid x_6=x_7=-x_8\}$.  We choose the basis of simple roots $\alpha_1 = \frac{1}{2}(\varepsilon_1+\varepsilon_8)-\frac{1}{2}(\varepsilon_2+\varepsilon_3+\varepsilon_4+\varepsilon_5+\varepsilon_6+\varepsilon_7), \alpha_2 = \varepsilon_1+\varepsilon_2, \alpha_3 = \varepsilon_2-\varepsilon_1, \alpha_4 = \varepsilon_3-\varepsilon_2, \alpha_5 = \varepsilon_4-\varepsilon_3, \alpha_6 = \varepsilon_5-\varepsilon_4$.  The positive roots are then
\[ \Phi^+ = \{ \pm \varepsilon_i+\varepsilon_j \mid 1 \leq i < j \leq 5\} \ \cup \ \left\{ \frac{1}{2}\left( \varepsilon_8-\varepsilon_7-\varepsilon_6+ \sum\limits_{i=1}^5 (-1)^{\nu(i)} \varepsilon_i \right) \ \middle| \  \sum\limits_{i=1}^5 \nu(i) \in 2\Z \right\}, \]
where $\nu(i) \in \{0,1\}$.
Table \ref{tab:E6} provides the expansion for each positive root in terms of $\Delta$. For roots in the second set, we record only the values $\nu(1), \dots, \nu(5)$ for brevity. 
 
 It is clear from the Dynkin diagram that $W/ \langle s_6 \rangle$ has type $D_5$.  We thus obtain the first 20 entries in Table \ref{tab:E6}, by reindexing the simple reflections $2 \mapsto 3, \ 3 \mapsto 4,\ 4 \mapsto 2$ from type $D_5$. Similarly, $W/ \langle s_1 \rangle$ has type $D_5$, and by reindexing the simple reflections $1 \mapsto 6, \ 2 \mapsto 5,\ 3 \mapsto 4, 4 \mapsto 2, \ 5 \mapsto 3$, ignoring all reflections in $W/\langle s_1, s_6 \rangle$ of type $D_4$ that already appear, we obtain the next 8 reflections in Table \ref{tab:E6}.  There is 1 additional reflection from the type $A_5$ system $W/ \langle s_2 \rangle$.  The 7 remaining reflections are distinct to type $E_6$, and are obtained by conjugating the 3 longest type $D_5$ reflections in $W/\langle s_1 \rangle$ by $s_1$, and then again by those subwords of $s_{243}$ which increase length. These reflections correspond to the 7 blue vertices in Figure \ref{fig:conj}.

\begin{table}[h]
\centering
  \begin{tabular}{ | c | c | c |}
    \hline
Root $\alpha \in \Phi^+$    & Expand $\alpha$ in $\Delta$ & Reflection $s_\alpha$ \\ \hline\hline
    $-\varepsilon_5+\varepsilon_6$    & $\alpha_7$ & $s_7$  \\ \hline 
       $-\varepsilon_4+\varepsilon_6$   & $\alpha_6+\alpha_7$ & $s_{767}$  \\ \hline 
        $-\varepsilon_3+\varepsilon_6$  & $\alpha_5+\alpha_6+\alpha_7$ & $s_{76567}$  \\ \hline
         $-\varepsilon_2+\varepsilon_6$ & $\alpha_4+\alpha_5+\alpha_6+\alpha_7$ & $s_{7654567}$  \\ \hline
                   $-\varepsilon_1+\varepsilon_6$& $\alpha_3+\alpha_4+\alpha_5+\alpha_6+\alpha_7$ & $s_{765434567}$  \\ \hline
         $\varepsilon_1+\varepsilon_6$ & $\alpha_2+\alpha_4+\alpha_5+\alpha_6+\alpha_7$ & $s_{765424567}$  \\ \hline

             $\varepsilon_2+\varepsilon_6$ & $\alpha_2+\alpha_3+\alpha_4+\alpha_5+\alpha_6+\alpha_7$ & $s_{76524342567}$  \\ \hline
                $\varepsilon_3+\varepsilon_6$   & $\alpha_2+\alpha_3+2\alpha_4+\alpha_5+\alpha_6+\alpha_7$ & $s_{4765243425674}$  \\ \hline
                 $\varepsilon_4+\varepsilon_6$  & $\alpha_2+\alpha_3+2\alpha_4+2\alpha_5+\alpha_6+\alpha_7$ & $s_{547652434256745}$  \\ \hline
             $\varepsilon_5+\varepsilon_6$   & $\alpha_2+\alpha_3+2\alpha_4+2\alpha_5+2\alpha_6+\alpha_7$ & $s_{65476524342567456}$  \\ \hline\hline
    111110 & $\alpha_1+\alpha_3+\alpha_4+\alpha_5+\alpha_6+\alpha_7$ & $s_{13456765431}$  \\ \hline\hline
        001110     & $\alpha_1+\alpha_2+\alpha_3+\alpha_4+\alpha_5+\alpha_6+\alpha_7$ & $s_{7165243425617}$ \\ \hline
       010110  & $\alpha_1+\alpha_2+\alpha_3+2\alpha_4+\alpha_5+\alpha_6+\alpha_7$ & $s_{714652434256417}$  \\ \hline
         011010   & $\alpha_1+\alpha_2+\alpha_3+2\alpha_4+2\alpha_5+\alpha_6+\alpha_7$ & $s_{71546524342564517}$  \\ \hline
           100110  & $\alpha_1+\alpha_2+2\alpha_3+2\alpha_4+\alpha_5+\alpha_6+\alpha_7$ & $s_{73146524342564137}$  \\ \hline
                   101010      & $\alpha_1+\alpha_2+2\alpha_3+2\alpha_4+2\alpha_5+\alpha_6+\alpha_7$ & $s_{7315465243425645137}$  \\ \hline
                    110010    & $\alpha_1+\alpha_2+2\alpha_3+3\alpha_4+2\alpha_5+\alpha_6+\alpha_7$ & $s_{743154652434256451347}$  \\ \hline
                   000010      & $\alpha_1+2\alpha_2+2\alpha_3+3\alpha_4+2\alpha_5+\alpha_6+\alpha_7$ & $s_{72431546524342564513427}$  \\ \hline
                     011100      & $\alpha_1+\alpha_2+\alpha_3+2\alpha_4+2\alpha_5+2\alpha_6+\alpha_7$ & $s_{6715465243425645176}$  \\ \hline
                   101100     & $\alpha_1+\alpha_2+2\alpha_3+2\alpha_4+2\alpha_5+2\alpha_6+\alpha_7$ & $s_{673154652434256451376}$  \\ \hline
                   110100     & $\alpha_1+\alpha_2+2\alpha_3+3\alpha_4+2\alpha_5+2\alpha_6+\alpha_7$ & $s_{67431546524342564513476}$  \\ \hline
                   000100      & $\alpha_1+2\alpha_2+2\alpha_3+3\alpha_4+2\alpha_5+2\alpha_6+\alpha_7$ & $s_{6724315465243425645134276}$  \\ \hline
                                 111000                 & $\alpha_1+\alpha_2+2\alpha_3+3\alpha_4+3\alpha_5+2\alpha_6+\alpha_7$ & $s_{5674315465243425645134765}$  \\ \hline
                                                                 001000          & $\alpha_1+2\alpha_2+2\alpha_3+3\alpha_4+3\alpha_5+2\alpha_6+\alpha_7$ & $s_{256743154652434256451347652}$  \\ \hline
                                                                                                               010000     & $\alpha_1+2\alpha_2+2\alpha_3+4\alpha_4+3\alpha_5+2\alpha_6+\alpha_7$ & $s_{42567431546524342564513476524}$  \\ \hline
   100000  & $\alpha_1+2\alpha_2+3\alpha_3+4\alpha_4+3\alpha_5+2\alpha_6+\alpha_7$ & $s_{3425674315465243425645134765243}$  \\ \hline
      $-\varepsilon_7+\varepsilon_8$  & $2\alpha_1+2\alpha_2+3\alpha_3+4\alpha_4+3\alpha_5+2\alpha_6+\alpha_7$ & $s_{134256743154652434256451347652431}$  \\ \hline
      \end{tabular} 
  \vskip 5pt
\caption{Reduced expressions for the 27 type $E_7$ reflections which are not type $E_6$.}\label{tab:E7}  
  \end{table}

    \subsection{Type $E_7$}

In type $E_7$, let $\varepsilon_i$ denote the projection of the standard basis vector $e_i \in \R^8$ onto onto the subspace $\{ (x_1, \dots, x_8) \mid x_7=-x_8\}$.  We choose the basis of simple roots $\alpha_1 = \frac{1}{2}(\varepsilon_1+\varepsilon_8)-\frac{1}{2}(\varepsilon_2+\varepsilon_3+\varepsilon_4+\varepsilon_5+\varepsilon_6+\varepsilon_7), \alpha_2 = \varepsilon_1+\varepsilon_2, \alpha_3 = \varepsilon_2-\varepsilon_1, \alpha_4 = \varepsilon_3-\varepsilon_2, \alpha_5 = \varepsilon_4-\varepsilon_3, \alpha_6 = \varepsilon_5-\varepsilon_4, \alpha_7 = \varepsilon_6-\varepsilon_5$.  The positive roots are then
\[ \Phi^+ = \{ \pm \varepsilon_i+\varepsilon_j \mid 1 \leq i < j \leq 6\}  \ \cup \ \{ \varepsilon_8-\varepsilon_7\}\ \cup \ \left\{ \frac{1}{2}\left( \varepsilon_8-\varepsilon_7+ \sum\limits_{i=1}^6 (-1)^{\nu(i)} \varepsilon_i \right) \ \middle| \  \sum\limits_{i=1}^6 \nu(i) \notin 2\Z \right\}, \]
where $\nu(i) \in \{0,1\}$.
Table \ref{tab:E7} provides the expansion for each positive root in terms of $\Delta$. For roots in the third set, we record only the values $\nu(1), \dots, \nu(6)$ for brevity.

 It is clear from both the choice of simple roots and the Dynkin diagram that $W/ \langle s_7 \rangle$ has type $E_6$, without any necessary relabeling.   We thus obtain the first 36 reflections in $E_7$ directly from Table \ref{tab:E6} for type $E_6$. Similarly, $W/ \langle s_1 \rangle$ has type $D_6$, and by reindexing the $D_6$ simple reflections $1 \mapsto 7, \ 2 \mapsto 6,\ 3 \mapsto 5, \ 5 \mapsto 2, 6 \mapsto 3$, ignoring all reflections in $W/\langle s_1, s_7 \rangle$ of type $D_5$ that have already been identified, we obtain the first 10 reflections in Table \ref{tab:E7}. There is 1 additional reflection not yet identified from the type $A_6$ system $W/ \langle s_2 \rangle$. The 16 remaining reflections are distinct to type $E_7$, and are obtained by conjugating the 7 type $E_6$ reflections in $W/\langle s_7 \rangle$ by $s_7$, then by $s_6$ provided that length increases, and finally by conjugating the second longest of the resulting reflections by $s_5, s_2, s_4, s_3, s_1$ in turn;  these reflections correspond to the 16 red vertices in Figure \ref{fig:conj}. The 27 type $E_7$ reflections additional to type $E_6$ appear in Table \ref{tab:E7}.

 \begin{table}[h]
\centering
  \begin{tabular}{ | c | c | c |}
    \hline
 Root $\alpha \in \Phi^+$     & Expand $\alpha$ in $\Delta$ & Reflection $s_\alpha$ \\ \hline\hline
   $-\varepsilon_6+\varepsilon_7$     & $\alpha_8$ & $s_8$  \\ \hline 
    $-\varepsilon_5+\varepsilon_7$      & $\alpha_7+\alpha_8$ & $s_{878}$  \\ \hline 
       $-\varepsilon_4+\varepsilon_7$   & $\alpha_6+\alpha_7+\alpha_8$ & $s_{87678}$  \\ \hline
      $-\varepsilon_3+\varepsilon_7$    & $\alpha_5+\alpha_6+\alpha_7+\alpha_8$ & $s_{8765678}$   \\ \hline
      $-\varepsilon_2+\varepsilon_7$    & $\alpha_4+\alpha_5+\alpha_6+\alpha_7+\alpha_8$ & $s_{876545678}$  \\ \hline
        
          $-\varepsilon_1+\varepsilon_7$ & $\alpha_3+\alpha_4+\alpha_5+\alpha_6+\alpha_7+\alpha_8$ &  $s_{87654345678}$ \\ \hline
           $\varepsilon_1+\varepsilon_7$  & $\alpha_2+\alpha_4+\alpha_5+\alpha_6+\alpha_7+\alpha_8$ & $s_{87654245678}$  \\ \hline
       $\varepsilon_2+\varepsilon_7$  & $\alpha_2+\alpha_3+\alpha_4+\alpha_5+\alpha_6+\alpha_7+\alpha_8$ & $s_{8765243425678}$  \\ \hline
       $\varepsilon_3+\varepsilon_7$  & $\alpha_2+\alpha_3+2\alpha_4+\alpha_5+\alpha_6+\alpha_7+\alpha_8$ & $s_{487652434256784}$  \\ \hline
         $\varepsilon_4+\varepsilon_7$ & $\alpha_2+\alpha_3+2\alpha_4+2\alpha_5+\alpha_6+\alpha_7+\alpha_8$ & $s_{54876524342567845}$  \\ \hline
           $\varepsilon_5+\varepsilon_7$ & $\alpha_2+\alpha_3+2\alpha_4+2\alpha_5+2\alpha_6+\alpha_7+\alpha_8$ & $s_{6548765243425678456}$  \\ \hline
           $\varepsilon_6+\varepsilon_7$ & $\alpha_2+\alpha_3+2\alpha_4+2\alpha_5+2\alpha_6+2\alpha_7+\alpha_8$ & $s_{765487652434256784567}$  \\ \hline\hline
   1111110   & $\alpha_1+\alpha_3+\alpha_4+\alpha_5+\alpha_6+\alpha_7+\alpha_8$ & $s_{1345678765431}$  \\ \hline
          \end{tabular} 
  \vskip 5pt
\caption{Reduced expressions for the 13 type $A_7$ and $D_7$ reflections in type $E_8$.}\label{tab:E8-1}  
  \end{table}

    \subsection{Type $E_8$}\label{sec:E8}
  
  In type $E_8$, we let $\varepsilon_i$ denote the standard basis vector in $\R^8$, for ease of comparison with types $E_6$ and $E_7$.  We choose the basis of simple roots $\alpha_1 = \frac{1}{2}(\varepsilon_1+\varepsilon_8)-\frac{1}{2}(\varepsilon_2+\varepsilon_3+\varepsilon_4+\varepsilon_5+\varepsilon_6+\varepsilon_7), \alpha_2 = \varepsilon_1+\varepsilon_2, \alpha_3 = \varepsilon_2-\varepsilon_1, \alpha_4 = \varepsilon_3-\varepsilon_2, \alpha_5 = \varepsilon_4-\varepsilon_3, \alpha_6 = \varepsilon_5-\varepsilon_4, \alpha_7 = \varepsilon_6-\varepsilon_5, \alpha_8 = \varepsilon_7-\varepsilon_6$. 
    The positive roots are then
\[ \Phi^+ = \{ \pm \varepsilon_i+\varepsilon_j \mid 1 \leq i < j \leq 8\} \ \cup \ \left\{ \frac{1}{2}\left( \varepsilon_8+ \sum\limits_{i=1}^7 (-1)^{\nu(i)} \varepsilon_i \right) \ \middle| \  \sum\limits_{i=1}^7 \nu(i) \in 2\Z \right\}, \]
where $\nu(i) \in \{0,1\}$. For roots in the second set, we record only the values $\nu(1), \dots, \nu(7)$ for brevity.
Tables \ref{tab:E8-1} and \ref{tab:E8-2} provide the expansion for each positive root in terms of $\Delta$.  Note that item (II) in the $E_8$ plate of \cite{Bour46} with these expansions contains 2 typos. Namely, writing $(c_1 \cdots c_8)$ for $\sum c_i\alpha_i \in \Phi^+$, the root $(12232211)$ is recorded twice, and an erroneous root $(11233321)$ appears. These roots are replaced by $(12232111)$ and $(11233221)$ in Table \ref{tab:E8-2} below.

 \begin{table}[ht]
\centering
  \begin{tabular}{ | c | c | c |}
    \hline
  Roots $\alpha \in \Phi^+$  & Expand $\alpha$ in $\Delta$ & $w$ for $s_\alpha = ws_\theta w^{-1}$  \\ \hline\hline
     0011110   & $\alpha_1+\alpha_2+\alpha_3+\alpha_4+\alpha_5+\alpha_6+\alpha_7+\alpha_8$  & $s_{87}$  \\ \hline 
       0101110 & $\alpha_1+\alpha_2+\alpha_3+2\alpha_4+\alpha_5+\alpha_6+\alpha_7+\alpha_8$& $s_{487}$  \\ \hline
       1001110 & $\alpha_1+\alpha_2+2\alpha_3+2\alpha_4+\alpha_5+\alpha_6+\alpha_7+\alpha_8$ & $s_{3487}$  \\ \hline
           1010110 & $\alpha_1+\alpha_2+2\alpha_3+2\alpha_4+2\alpha_5+\alpha_6+\alpha_7+\alpha_8$ & $s_{53487}$  \\ \hline 
                 1100110  & $\alpha_1+\alpha_2+2\alpha_3+3\alpha_4+2\alpha_5+\alpha_6+\alpha_7+\alpha_8$ & $s_{453487}$  \\ \hline 
                      0000110  & $\alpha_1+2\alpha_2+2\alpha_3+3\alpha_4+2\alpha_5+\alpha_6+\alpha_7+\alpha_8$ & $s_{2453487}$  \\ \hline 
                       0001010     & $\alpha_1+2\alpha_2+2\alpha_3+3\alpha_4+2\alpha_5+2\alpha_6+\alpha_7+\alpha_8$ & $s_{62453487}$  \\ \hline 
        0110110  &$\alpha_1+\alpha_2+\alpha_3+2\alpha_4+2\alpha_5+\alpha_6+\alpha_7+\alpha_8$ & $s_{5487}$  \\ \hline
     0111010  & $\alpha_1+\alpha_2+\alpha_3+2\alpha_4+2\alpha_5+2\alpha_6+\alpha_7+\alpha_8$ & $s_{65487}$  \\ \hline 
        1011010  & $\alpha_1+\alpha_2+2\alpha_3+2\alpha_4+2\alpha_5+2\alpha_6+\alpha_7+\alpha_8$  & $s_{365487}$  \\ \hline 
         1101010    &  $\alpha_1+\alpha_2+2\alpha_3+3\alpha_4+2\alpha_5+2\alpha_6+\alpha_7+\alpha_8$ & $s_{4365487}$  \\ \hline 
               1110010  &$\alpha_1+\alpha_2+2\alpha_3+3\alpha_4+3\alpha_5+2\alpha_6+\alpha_7+\alpha_8$ & $s_{54365487}$  \\ \hline 
                0010010    & $\alpha_1+2\alpha_2+2\alpha_3+3\alpha_4+3\alpha_5+2\alpha_6+\alpha_7+\alpha_8$ & $s_{254365487}$  \\ \hline 
                         0100010      &$\alpha_1+2\alpha_2+2\alpha_3+4\alpha_4+3\alpha_5+2\alpha_6+\alpha_7+\alpha_8$ & $s_{4254365487}$  \\ \hline 
                                      1000010    & $\alpha_1+2\alpha_2+3\alpha_3+4\alpha_4+3\alpha_5+2\alpha_6+\alpha_7+\alpha_8$ & $s_{34254365487}$  \\ \hline 
                                       $-\varepsilon_6+\varepsilon_8$              &$2\alpha_1+2\alpha_2+3\alpha_3+4\alpha_4+3\alpha_5+2\alpha_6+\alpha_7+\alpha_8$ & $s_{134254365487}$  \\ \hline 
       0111100     & $\alpha_1+\alpha_2+\alpha_3+2\alpha_4+2\alpha_5+2\alpha_6+2\alpha_7+\alpha_8$ & $s_{765487}$  \\ \hline     
         1011100  & $\alpha_1+\alpha_2+2\alpha_3+2\alpha_4+2\alpha_5+2\alpha_6+2\alpha_7+\alpha_8$  & $s_{3765487}$  \\ \hline
           1101100  & $\alpha_1+\alpha_2+2\alpha_3+3\alpha_4+2\alpha_5+2\alpha_6+2\alpha_7+\alpha_8$ & $s_{43765487}$  \\ \hline
       0001100  & $\alpha_1+2\alpha_2+2\alpha_3+3\alpha_4+2\alpha_5+2\alpha_6+2\alpha_7+\alpha_8$ & $s_{243765487}$  \\ \hline
         0010100   & $\alpha_1+2\alpha_2+2\alpha_3+3\alpha_4+3\alpha_5+2\alpha_6+2\alpha_7+\alpha_8$& $s_{5243765487}$  \\ \hline
         0100100           & $\alpha_1+2\alpha_2+2\alpha_3+4\alpha_4+3\alpha_5+2\alpha_6+2\alpha_7+\alpha_8$ & $s_{45243765487}$  \\ \hline
              1000100           & $\alpha_1+2\alpha_2+3\alpha_3+4\alpha_4+3\alpha_5+2\alpha_6+2\alpha_7+\alpha_8$ & $s_{345243765487}$  \\ \hline
                  $-\varepsilon_5+\varepsilon_8$             & $2\alpha_1+2\alpha_2+3\alpha_3+4\alpha_4+3\alpha_5+2\alpha_6+2\alpha_7+\alpha_8$ & $s_{1345243765487}$  \\ \hline
                1110100                & $\alpha_1+\alpha_2+2\alpha_3+3\alpha_4+3\alpha_5+2\alpha_6+2\alpha_7+\alpha_8$ & $s_{543765487}$  \\ \hline
       1111000     & $\alpha_1+\alpha_2+2\alpha_3+3\alpha_4+3\alpha_5+3\alpha_6+2\alpha_7+\alpha_8$ & $s_{6543765487}$  \\ \hline
        0011000         & $\alpha_1+2\alpha_2+2\alpha_3+3\alpha_4+3\alpha_5+3\alpha_6+2\alpha_7+\alpha_8$ & $s_{26543765487}$  \\ \hline
         0101000     & $\alpha_1+2\alpha_2+2\alpha_3+4\alpha_4+3\alpha_5+3\alpha_6+2\alpha_7+\alpha_8$& $s_{426543765487}$  \\ \hline
        1001000    & $\alpha_1+2\alpha_2+3\alpha_3+4\alpha_4+3\alpha_5+3\alpha_6+2\alpha_7+\alpha_8$ & $s_{3426543765487}$  \\ \hline
           $-\varepsilon_4+\varepsilon_8$    & $2\alpha_1+2\alpha_2+3\alpha_3+4\alpha_4+3\alpha_5+3\alpha_6+2\alpha_7+\alpha_8$ & $s_{13426543765487}$  \\ \hline
               $-\varepsilon_3+\varepsilon_8$      & $2\alpha_1+2\alpha_2+3\alpha_3+4\alpha_4+4\alpha_5+3\alpha_6+2\alpha_7+\alpha_8$ & $s_{513426543765487}$  \\ \hline
                 $-\varepsilon_2+\varepsilon_8$       & $2\alpha_1+2\alpha_2+3\alpha_3+5\alpha_4+4\alpha_5+3\alpha_6+2\alpha_7+\alpha_8$ & $s_{4513426543765487}$  \\ \hline
                   $-\varepsilon_1+\varepsilon_8$      & $2\alpha_1+2\alpha_2+4\alpha_3+5\alpha_4+4\alpha_5+3\alpha_6+2\alpha_7+\alpha_8$ & $s_{34513426543765487}$  \\ \hline
                    0110000         &$\alpha_1+2\alpha_2+2\alpha_3+4\alpha_4+4\alpha_5+3\alpha_6+2\alpha_7+\alpha_8$ & $s_{5426543765487}$  \\ \hline
        1010000   & $\alpha_1+2\alpha_2+3\alpha_3+4\alpha_4+4\alpha_5+3\alpha_6+2\alpha_7+\alpha_8$ & $s_{35426543765487}$  \\ \hline
        1100000      & $\alpha_1+2\alpha_2+3\alpha_3+5\alpha_4+4\alpha_5+3\alpha_6+2\alpha_7+\alpha_8$ & $s_{435426543765487}$  \\ \hline
          0000000   & $\alpha_1+3\alpha_2+3\alpha_3+5\alpha_4+4\alpha_5+3\alpha_6+2\alpha_7+\alpha_8$ & $s_{2435426543765487}$  \\ \hline
              $\varepsilon_1+\varepsilon_8$    & $2\alpha_1+3\alpha_2+3\alpha_3+5\alpha_4+4\alpha_5+3\alpha_6+2\alpha_7+\alpha_8$ & $s_{12435426543765487}$  \\ \hline
                  $\varepsilon_2+\varepsilon_8$   & $2\alpha_1+3\alpha_2+4\alpha_3+5\alpha_4+4\alpha_5+3\alpha_6+2\alpha_7+\alpha_8$ & $s_{312435426543765487}$  \\ \hline
                     $\varepsilon_3+\varepsilon_8$     & $2\alpha_1+3\alpha_2+4\alpha_3+6\alpha_4+4\alpha_5+3\alpha_6+2\alpha_7+\alpha_8$ & $s_{4312435426543765487}$  \\ \hline
                    $\varepsilon_4+\varepsilon_8$      & $2\alpha_1+3\alpha_2+4\alpha_3+6\alpha_4+5\alpha_5+3\alpha_6+2\alpha_7+\alpha_8$ & $s_{54312435426543765487}$  \\ \hline
                      $\varepsilon_5+\varepsilon_8$    & $2\alpha_1+3\alpha_2+4\alpha_3+6\alpha_4+5\alpha_5+4\alpha_6+2\alpha_7+\alpha_8$  & $s_{654312435426543765487}$  \\ \hline
 $\varepsilon_6+\varepsilon_8$ & $2\alpha_1+3\alpha_2+4\alpha_3+6\alpha_4+5\alpha_5+4\alpha_6+3\alpha_7+\alpha_8$ & $s_{7654312435426543765487}$  \\ \hline
  $\varepsilon_7+\varepsilon_8$ & $2\alpha_1+3\alpha_2+4\alpha_3+6\alpha_4+5\alpha_5+4\alpha_6+3\alpha_7+2\alpha_8$ & $s_{87654312435426543765487}$  \\ \hline 
          \end{tabular} 
  \vskip 5pt
\caption{Reduced expressions for 44 reflections distinct to type $E_8$, with $s_{\theta} = s_{16524342561}$.}\label{tab:E8-2}  
  \end{table}

 It is clear from both the choice of simple roots and the Dynkin diagram that $W/ \langle s_8 \rangle$ has type $E_7$, without any necessary relabeling.   We thus obtain the first 63 reflections in $E_8$ directly from Tables \ref{tab:E6} and \ref{tab:E7} displaying reduced words for the reflections in types $E_6$ and $E_7$, respectively.  In addition, $W/ \langle s_1 \rangle$ has type $D_7$, and by reindexing the $D_7$ simple reflections $1 \mapsto 8, \ 2 \mapsto 7,\ 3 \mapsto 6,\ 4 \mapsto 5, \ 5 \mapsto 4, 6 \mapsto 2,\ 7 \mapsto 3$, ignoring all reflections in $W/\langle s_1, s_8 \rangle$ of type $D_6$ that have already been identified, we obtain the first 12 reflections in Table \ref{tab:E8-1}. There is one additional reflection from the type $A_7$ system $W/ \langle s_2 \rangle$, which has both $s_1$ and $s_8$ in its support, shown as the last entry in Table \ref{tab:E8-1}.

  The 44 remaining reflections are distinct to type $E_8$, and can be obtained by conjugating the 16 type $E_7$ reflections in $W/\langle s_8 \rangle$ by $s_8$, and then performing additional conjugations as depicted by Figure \ref{fig:conj}.  As explained in Example \ref{ex:introE8}, the vertices of the graph in Figure \ref{fig:conj} correspond to the 67 reflections of types $E_6, E_7,$ and $E_8$, colored blue, red, and green respectively, with length increasing going upward. The type $E_6$ reflection of minimal length is  $s_{\theta} = s_{16524342561}$, obtained from the bottom portion of Table \ref{tab:E6}.  The edge labels in Figure \ref{fig:conj} indicate the simple reflection with which to conjugate to obtain the reflection indexing the adjacent vertex. For example, conjugating $s_\theta$ by $w =s_{24354}$ reading up the blue edges along the righthand side gives one reduced expression for the longest reflection in type $E_6$; compare the final entry of Table \ref{tab:E6}.  Similarly, the red vertices in Figure \ref{fig:conj} correspond to the 17 reflections of type $E_7$, and the green vertices correspond to the 44 remaining reflections which are distinct to type $E_8$. 
 Starting from $s_{87}s_\theta s_{78}$ at the bottom of the green vertices and conjugating by the labels on the green edges, we obtain the 44 reduced expressions recorded in Table \ref{tab:E8-2} on the next page.

  This concludes our proof of Theorem \ref{thm:main} in the exceptional types. \qed


\bibliographystyle{alphanum}
\bibliography{references}

\end{document}